\documentclass{article}
\usepackage{multirow}
\usepackage{tabularx}
\usepackage{PRIMEarxiv}
\usepackage{times}
\usepackage{fullpage,comment}
\usepackage{amsfonts}
\usepackage{color}
\usepackage[sort]{cite}
\usepackage{graphicx}
\usepackage[utf8]{inputenc}
\usepackage[T1]{fontenc}
\usepackage[mathlines]{lineno}
\usepackage{amssymb,amsthm,amsmath,bm}
\usepackage{hyperref} 
\usepackage{algorithmic}
\usepackage{paralist} 
\usepackage{booktabs} 
\usepackage[algo2e,ruled,noend,linesnumbered]{algorithm2e}

\usepackage{authblk} 

\def\d{{\, \rm d}}
\def\u{\mathbf{u}}

\def\y{\mathbf{y}}
\def\F{\mathbf{f}}

\def\R{\mathbf{R}}
\def\X{\mathbf{X}}
\def\Y{\mathbf{Y}}

\def\btheta{\bm{\theta}}
\def\bmu{\bm{\mu}}

\def\btheta{\bm{\theta}}

\definecolor{orange}{RGB}{255,127,0}

\definecolor{purple}{RGB}{76,0,153}

\def\peq{p_{\textrm{eq}}}

\def\cQ{\mathcal{Q}}
\def\obs{\textrm{obs}}

\usepackage{lineno}

\def\u{\mathbf{u}}
\def\A{\mathbf{A}}

\def\a{\mathbf{a}}

\def\B{\mathbf{B}}
\def\b{\mathbf{b}}
\def\c{\mathbf{c}}

\def\D{\mathbf{D}}

\def\F{\mathbf{F}}

\def\G{\mathbf{G}}

\def\M{\mathbf{M}}

\def\R{\mathbf{R}}
\def\S{\mathbf{S}}

\def\u{\mathbf{u}}

\def\w{\mathbf{w}}

\def\X{\mathbf{X}}
\def\y{\mathbf{y}}
\def\Y{\mathbf{Y}}

\def\z{\mathbf{z}}

\def\cB{\mathcal{B}}
\def\cL{\mathcal{L}}

\def\cQ{\mathcal{Q}}

\def\bxi{\bm{\xi}}
\def\btheta{\bm{\theta}}

\def\peq{p_{\textrm{eq}}}

\def\bmu{\bm{\mu}}

\title{Conditional Gaussian Nonlinear System: a Fast Preconditioner and a Cheap Surrogate Model For  Complex Nonlinear Systems} 

\author[1]{Nan Chen}
\author[1]{Yingda Li}
\author[2]{Honghu Liu}
\affil[1]{Department of Mathematics, University of Wisconsin-Madison, Madison, WI 53705, USA (chennan@math.wisc.edu; yli678@wisc.edu)}
\affil[2]{Department of Mathematics, Virginia Tech, Blacksburg, VA 24061, USA (hhliu@vt.edu)}

\date{\today}
\begin{document}\maketitle
\begin{abstract}
Developing suitable approximate models for analyzing and simulating complex nonlinear systems is practically important. This paper aims at exploring the skill of a rich class of nonlinear stochastic models, known as the conditional Gaussian nonlinear system (CGNS), as both a cheap surrogate model and a fast preconditioner for facilitating many computationally challenging tasks. The CGNS preserves the underlying physics to a large extent and can reproduce intermittency, extreme events and other non-Gaussian features in many complex systems arising from practical applications. Three interrelated topics are studied. First, the closed analytic formulae of solving the conditional statistics provide an efficient and accurate data assimilation scheme. It is shown that the data assimilation skill of a suitable CGNS approximate forecast model outweighs that by applying an ensemble method even to the perfect model with strong nonlinearity, where the latter suffers from filter divergence. Second, the CGNS allows the development of a fast algorithm for simultaneously estimating the parameters and the unobserved variables with uncertainty quantification in the presence of only partial observations. Utilizing an appropriate CGNS as a preconditioner significantly reduces the computational cost in accurately estimating the parameters in the original complex system. Finally, the CGNS advances rapid and statistically accurate algorithms for computing the probability density function and sampling the trajectories of the unobserved state variables.  These fast algorithms facilitate the development of an efficient and accurate data-driven method for predicting the linear response of the original system with respect to parameter perturbations based on a suitable CGNS preconditioner.
\end{abstract}

\section{Introduction}

Complex nonlinear systems are ubiquitous in many areas, including geophysics, climate science, engineering, neuroscience and material science~\cite{vallis2017atmospheric, strogatz2018nonlinear, wilcox1988multiscale, sheard2009principles,ghil2012topics}. Mathematical modeling plays an important role in characterizing and discovering the underlying physics of these complex systems~\cite{edwards2011history, majda2012challenges}.
With suitable mathematical models in hand, effective parameter inference, state estimation and data assimilation  become fundamental tasks that serve as the prerequisites for analyzing these systems~\cite{asch2016data, kalnay2003atmospheric, majda2012filtering, law2015data, ghil1991data}. Accurate forecast of future states and successful prediction of the system response due to external perturbations are also central topics that have many practical implications~\cite{lucarini2017predicting, ragone2016new, majda2010low, toth1997ensemble, leutbecher2008ensemble}.

However, there exist quite a few mathematical and computational challenges in analyzing and simulating complex nonlinear systems. First, the intrinsic nonlinearity in these complex nonlinear systems often triggers strongly chaotic or turbulent behavior~\cite{salmon1998lectures, dijkstra2013nonlinear, palmer1993nonlinear}. As a consequence, intermittency, extreme events and non-Gaussian probability density functions (PDFs) are some of the typical features in these systems~\cite{farazmand2019extreme, trenberth2015attribution, moffatt2021extreme, majda2003introduction, manneville1979intermittency}, which  impede the use of many traditional mathematical tools to analyze the model properties. Second, due to the nonlinear interactions between state variables across different scales, many of these complex nonlinear systems are high dimensional and have multiscale spatiotemporal structures~\cite{wilcox1988multiscale, majda2016introduction, tao2009multiscale, majda2014new}. Therefore, developing new efficient numerical algorithms to accelerate the computational efficiency becomes essential. Particularly, enhancing the computational efficiency by reducing the complexity of these systems via effective stochastic parameterizations is practically important and is widely used in for example climate sciences. Third, it is often the case in practice that only partial observations of the state variables are available~\cite{kalnay2003atmospheric, lau2011intraseasonal}, which result in additional difficulties for model calibration, state estimation and prediction where systematic uncertainty quantification needs to be addressed~\cite{curry2011climate, delsole2004predictability, edwards1999global, majda2012lessons, majda2018model}.

Since many complex dynamical systems are too expensive to be handled directly, it is of practical importance to develop suitable approximate models, which capture certain features of nature and are easier to deal with. Exploiting systematic reduced order modeling strategies and effective (stochastic) parameterizations is often a pre-requisite for the development of approximate models. Linear regression models are arguably the simplest class of approximate models~\cite{yan2009linear, freedman2009statistical}, which can already provide certain skill for short-term forecasts although they usually suffer in characterizing the underlying nonlinear physics. Physics-constrained regression models are a set of nonlinear approximate models~\cite{majda2012physics, harlim2014ensemble,kondrashov2015data}, which take into account the energy conserving nonlinear interactions in the model development that guarantees the well-posednss of long-term behavior of the system. Another commonly used approach to developing approximate models is to project the starting complex nonlinear system to the leading a few energetic modes in light of the Galerkin proper orthogonal decomposition methods~\cite{HLB96} or other empirical basis functions such as the principal interaction patterns \cite{hasselmann1988pips,kwasniok1996reduction} and the dynamic mode decomposition \cite{rowley2009spectral,schmid2010dynamic}. With a careful design of the closure terms to compensate the truncation error, these reduced order models are skillful in resolving certain problems in fluids and turbulence~\cite{carlberg2013gnat, noack2011reduced, taira2020modal, xie2018data}.

Meanwhile, many data-driven approximate modeling strategies have recently been developed~\cite{ahmed2021closures, chekroun2017data, chekroun2021stochastic, lin2021data, mou2021data, peherstorfer2015dynamic, hijazi2020data, smarra2018data}. One of them is the sparse identification of nonlinear dynamical systems (SINDy)~\cite{brunton2016discovering}, which leads to nonlinear regression models with parsimonious structures via sparse regression and compressed sensing. Many other approximate modeling approaches have also been designed for specific scientific purposes. For example, the past noise forecasting method~\cite{chekroun2011predicting} was developed as a data-driven forecast model for stochastic climate processes that exhibit low-frequency variability. Reduced-space Gaussian process regression forecast~\cite{wan2017reduced} was designed for data-driven probabilistic forecast of chaotic dynamical systems. Small-scale parameterization based on a data-informed optimal homotopic deformation of invariant manifolds was developed to design low-dimensional models for both deterministic chaotic systems and stochastic systems \cite{CLW15_vol2,CLM20}. Physically consistent data-driven weather forecasting techniques were proposed and applied to operational models~\cite{chattopadhyay2020deep, chattopadhyay2021towards}. Recently, a strong link between the stochastic parameterization approach based on perturbation expansions of the Koopman operator \cite{wouters2013multi} and the data-driven empirical model reduction (EMR) methodology \cite{kravtsov2005multilevel} was established in~\cite{santos2021reduced}. In addition, machine learning methods nowadays have been extensively incorporated into the reduced order models to further improve the approximation and forecast skill~\cite{san2018extreme, chattopadhyay2020superparameterization, chen2021bamcafe, pawar2020data, moosavi2015efficient}.

The objective of this paper is to explore the skill of a rich class of nonlinear stochastic models, known as the ``conditional Gaussian nonlinear system'' (CGNS)~\cite{chen2018conditional}, as approximate models for  complex nonlinear systems. The CGNS includes many physics-constrained nonlinear stochastic models (e.g., the stochastic versions of various Lorenz models, low-order models of Charney-DeVore flows, and a paradigm model for topographic mean flow interaction), quite a few stochastically coupled reaction-diffusion models in neuroscience and ecology (e.g., stochastically coupled FitzHugh-Nagumo (FHN) models and stochastically coupled SIR epidemic models), and several large-scale dynamical models in engineering and geophysical flows (e.g., the Boussinesq equations with noise and stochastically forced rotating shallow water equation). See the  article~\cite{chen2018conditional} for a gallery of examples of the CGNS. The CGNS has also been applied to modeling and forecasting several important climate phenomena, such as the Madden-Julian oscillation and the monsoon~\cite{chen2014predicting, chen2018predicting}, and has been utilized for Lagrangian data assimilation~\cite{chen2014information}. Yet, most of the previous work focused on perfect model scenarios, where utilizing the CGNS as an approximate model has not been systematically studied.

The CGNS has several unique features that allows it to be distinct from many existing approximate modeling strategies. First, the CGNS aims at preserving the underlying physical mechanism to the greatest extent. Specifically, the nonlinearity involving the large-scale or slow variables are by design retained, which includes not only the self-interactions among the large-scale variables but also the cross-scale interactions between large- and small-scale variables, while suitable approximations are imposed primarily on the nonlinear self-interactions between small-scale, fast-varying or unresolved state variables via effective stochastic parameterizations. This is fundamentally different from many purely data-driven nonlinear regression models, which may miss certain crucial underlying physics of the original complex nonlinear systems. Second, the stochastic parameterizations of the self-interactions between small-scale variables lead to an important feature of the CGNS. That is, the distribution of the small-scale variables conditioned on the large-scale ones is Gaussian. One remarkable consequence is that the associated conditional Gaussian distribution can be calculated using closed analytic formulae~\cite{liptser2013statistics}, which considerably facilitate the mathematical analysis and numerical simulations of the CGNS. In fact, the closed analytic formulae of the conditional Gaussian distributions allow the development of efficient and statistically accurate algorithms for parameter estimation, data assimilation and ensemble forecast in light of only partial observations. Note that, despite the conditional Gaussianity, the joint and marginal distributions of the CGNS remain highly non-Gaussian. Thus, the intermittency, extreme events and turbulent features can all be preserved in a suitably designed CGNS. Third, the CGNS is also adaptable to many data-driven scenarios.  Physics constraints, localizations and sparse identification together with many other mathematical and computational strategies can be possibly incorporated into the CGNS. Finally, information theory~\cite{kleeman2011information, majda2005information} can be applied to quantify the uncertainty and the statistical error of the CGNS in approximating the original complex nonlinear systems.

\begin{figure}[tbh!]
\centering
\includegraphics[width=0.95\textwidth]{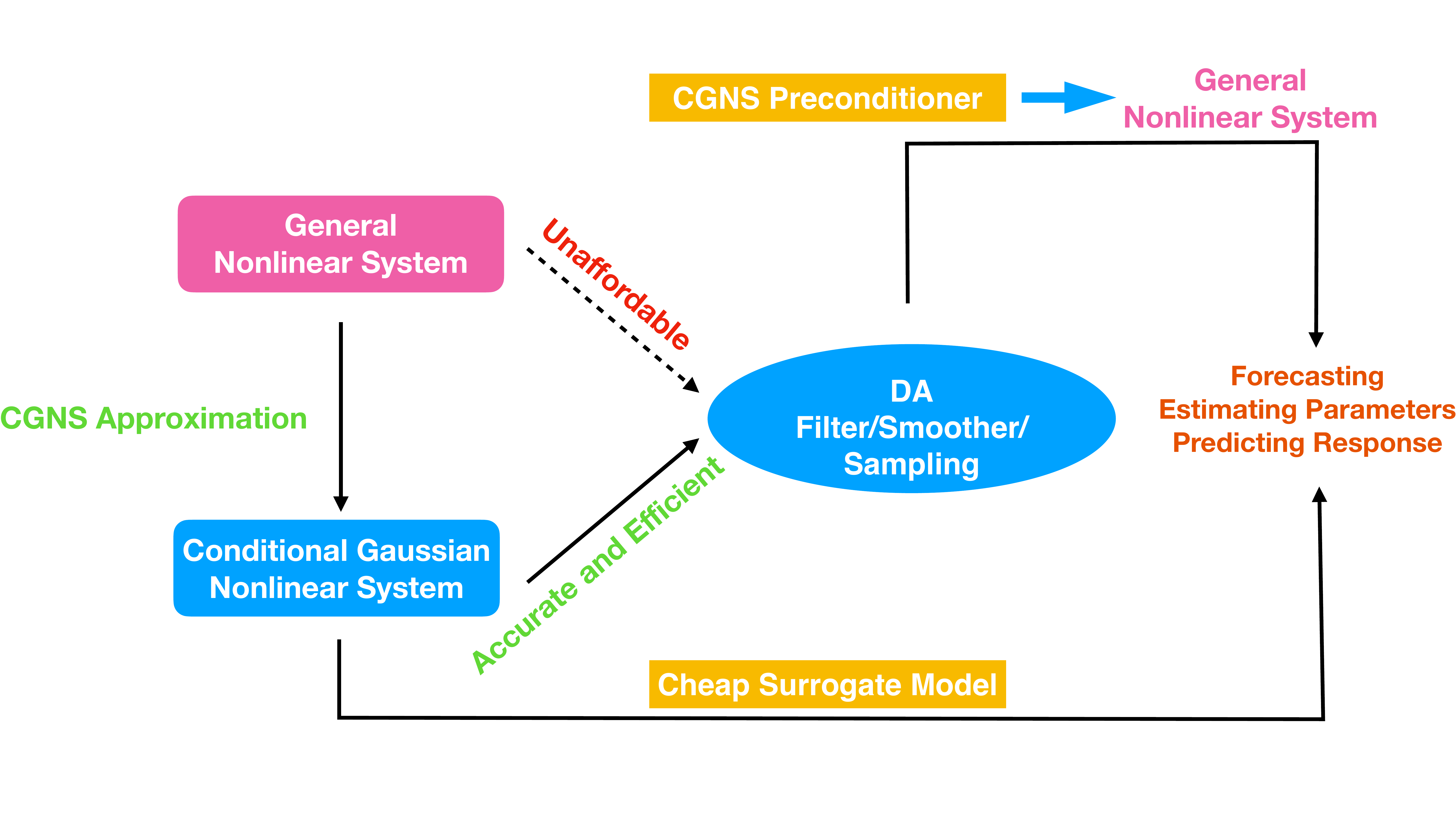}
  \caption{A schematic illustration of utilizing the CGNS as a fast preconditioner and a cheap surrogate model for general nonlinear system.}\label{fig: CGNS_illustration}
\end{figure}

The specific goal of this paper is two-fold. First, by taking advantage of its analytically solvable properties, the CGNS can be served as a fast preconditioner for facilitating many computationally challenging tasks associated with the original complex nonlinear system. Important applications include estimating the parameters of the original system in the presence of only partial observations and recovering the non-Gaussian PDFs as a crucial intermediate step for computing the model sensitivity and response.
Second, the CGNS is exploited as a surrogate model by exploiting systematic reduced order modeling strategies and suitable stochastic parameterizations, aiming at spending a much lower computational cost to create comparably accurate results as those obtained from the original complex system. This includes for example the state estimation of unobserved variables and the statistical forecast. Figure~\ref{fig: CGNS_illustration} shows a schematic illustration of utilizing the CGNS as a fast preconditioner and a cheap surrogate model for general nonlinear system.

The rest of the paper is organized as follows. The general mathematical framework of CGNS and its analytic properties are described in Section~\ref{Sec:CGNS}. Several systematic strategies for the development of the CGNS are included in Section~\ref{Sec:Model_Development}. Sections~\ref{Sec:DA_Forecast} --~\ref{Sec:Response} consist of three important tasks in complex nonlinear systems, showing the roles of the CGNS both as a surrogate model and a preconditioner for the original complex system. Specifically, Section~\ref{Sec:DA_Forecast} focuses on the data assimilation and ensemble forecast, Section~\ref{Sec:Parameter_Estimation} aims at efficient parameter estimation, and Section~\ref{Sec:Response} illustrates the use of CGNS in facilitating the study of model sensitivity and response theory. The paper is concluded in~Section~\ref{Conclusion}.

\section{General Mathematical Framework of the CGNS}\label{Sec:CGNS}

Let us start with the general formulation of the turbulent dynamical systems motivated from fluid and geophysical applications~\cite{vallis2017atmospheric, salmon1998lectures, kalnay2003atmospheric, majda2016introduction},
\begin{equation}\label{eq:abs_formu}
\frac{\mathrm{d}\u}{\mathrm{d}t}=\left(L+D\right)\u+B\left(\u,\u\right)+\F\left(t\right)+\boldsymbol{\sigma}\left(\u, t\right)\dot{\mathbf{W}}\left(t\right),
\end{equation}
where the state variable $\u\in\mathbb{C}^{N}$ is in a high dimensional phase space.  In~\eqref{eq:abs_formu}, the first two components, $\left(L+D\right)\u$,
represent linear dispersion and dissipation effects, where $L^{*}=-L$ is a skew-symmetric operator (with $\cdot^*$ being the complex conjugate transpose); and $D$ is a negative-definite matrix. The nonlinear effect is introduced through an energy-conserving quadratic form, $B\left(\u,\u\right)$. Besides, the system is subject to external forcing effects that are decomposed into a deterministic component, $\F\left(t\right)$, and a stochastic component represented by a Gaussian random process, $\boldsymbol{\sigma}\left(\u, t\right)\dot{\mathbf{W}}\left(t\right)$, where $\boldsymbol{\sigma}\in \mathbb{C}^{N\times K}$ is the noise matrix and $\dot{\mathbf{W}}\in\mathbb{C}^K$ is the white noise. The two components $\left(L+D\right)\u+B\left(\u,\u\right)+\F\left(t\right)$ and $\boldsymbol{\sigma}\left(\u, t\right)$ on the right hand side of~\eqref{eq:abs_formu} are also known as the drift part and the diffusion coefficients, respectively.

\subsection{The CGNS}
Despite being highly nonlinear and possessing strongly non-Gaussian statistics in both the marginal and joint distributions of the state $\u$, many complex nonlinear dynamical systems~\eqref{eq:abs_formu} have or can be  approximated by the following nonlinear system with conditional Gaussian structures. The general mathematical framework of the CGNS is  as follows~\cite{liptser2013statistics, chen2016filtering, chen2018conditional},
\begin{subequations}\label{CGNS}
\begin{align}
  \frac{\d\mathbf{X}}{\d t} &= \Big[\mathbf{A}_\mathbf{0}(\mathbf{X},t) + \mathbf{A}_\mathbf{1}(\mathbf{X},t) \mathbf{Y}(t)\Big]  + \mathbf{B}_\mathbf{1}(\mathbf{X},t)\dot{\mathbf{W}}_\mathbf{1}(t),\label{CGNS_X}\\
  \frac{\d\mathbf{Y}}{\d t} &= \Big[\mathbf{a}_\mathbf{0}(\mathbf{X},t) + \mathbf{a}_\mathbf{1}(\mathbf{X},t) \mathbf{Y}(t)\Big]  + \mathbf{b}_\mathbf{2}(\mathbf{X},t)\dot{\mathbf{W}}_\mathbf{2}(t),\label{CGNS_Y}
\end{align}
\end{subequations}
where the original model state $\u$ is decomposed into multi-dimensional state variables $\mathbf{X}\in\mathbb{C}^{N_1}$ and $\mathbf{Y}\in\mathbb{C}^{N_2}$, with $N_1+N_2=N$. In~\eqref{CGNS}, $\mathbf{A}_0, \mathbf{a}_0, \mathbf{A}_1, \mathbf{a}_1, \mathbf{B}_1$ and $\mathbf{b}_2$ are vectors or matrices that can depend nonlinearly on the state variables $\mathbf{X}$ and time $t$  while $\dot{\mathbf{W}}_\mathbf{1}$ and $\dot{\mathbf{W}}_\mathbf{2}$ are independent white noise sources that can have different dimensions from $\mathbf{X}$ and $\mathbf{Y}$.

The name ``conditional Gaussian'' comes from the fact that once a time series of $\mathbf{X}(s)$ for $s\leq t$ is given, then the conditional distribution $p(\mathbf{Y}(t)|\mathbf{X}(s\leq t))$ is Gaussian. This can be seen by noticing that, with a given $\mathbf{X}$, the process of $\mathbf{Y}$ is linear (with respect to the variable $\mathbf{Y}$ itself since $\mathbf{X}$ has been given) with Gaussian white noises. It is worthwhile to highlight that, from the general form of the complex turbulent system~\eqref{eq:abs_formu} to the CGNS~\eqref{CGNS}, the nonlinear interactions between $\mathbf{X}$ itself and those between $\mathbf{X}$ and $\mathbf{Y}$ in~\eqref{eq:abs_formu} are both completely retained. The only  simplification in the CGNS is to approximate the  nonlinear self-interactions between $\mathbf{Y}$ by a combination of nonlinear functions of $\mathbf{X}$, conditional linear functions of $\mathbf{Y}$ and effective stochastic noises. Nevertheless, if $\mathbf{Y}$ represents small-scale or fast variables, then such a manipulation  is expected to be an effective approximation that preserves the underlying physics to a large extent.

It should be noted that the CGNS in~\eqref{CGNS} is still highly nonlinear due to the nonlinearity in $\mathbf{A}_0, \mathbf{a}_0, \mathbf{A}_1$ and $\mathbf{a}_1$ as well as the nonlinear coupling between the latter two with $\mathbf{Y}$. Such nonlinearities preserve the non-Gaussian statistics in~\eqref{eq:abs_formu} and allow to reproduce many observed features in nature such as extreme events with the more tractable conditional Gaussian structure.
A gallery of examples of the CGNS, including many physics-constrained nonlinear stochastic models, quite a few stochastically coupled reaction-diffusion models in neuroscience and ecology, and some large-scale dynamical models in engineering and geophysical flows can be found in~\cite{chen2018conditional}.

Despite being highly nonlinear and non-Gaussian, one of the important features of the CGNS~\eqref{CGNS} is that the conditional distribution of $\mathbf{Y}$ given one realization of the time series $\mathbf{X}$ can be solved via closed analytic formulae. Such a unique analytic property significantly facilitates the analysis and calculations of state estimation, data assimilation and forecast. This feature also makes the CGNS to be quite different from general nonlinear or non-Gaussian systems. For the latter, particle methods have to be applied for state estimation \cite{carpenter1999improved, hol2006resampling, nummiaro2003adaptive}, in which many empirical tunings have to be involved to mitigate the numerical sampling errors.

\subsection{Closed analytic formulae for computing the conditional statistics and data assimilation}\label{Subsec:Pathwise}
\subsubsection{Nonlinear filter}
For the CGNS~\eqref{CGNS}, given one realization of the time series $\mathbf{X}(s)$ for $s\in[0,t]$, the conditional distribution
\begin{equation}\label{CGNS_PDF}
    p(\mathbf{Y}(t)|\mathbf{X}(s),s\leq t) \sim \mathcal{N}(\boldsymbol\mu_{\mathbf{f}}(t),\mathbf{R}_{\mathbf{f}}(t))
\end{equation}
becomes Gaussian, where the conditional mean $\boldsymbol\mu_{\mathbf{f}}$ and the conditional covariance $\mathbf{R}_{\mathbf{f}}$ are given by the following explicit formulae~\cite{liptser2013statistics}
\begin{subequations}\label{CGNS_Stat}
\begin{align}
  \frac{\d \boldsymbol{\mu}_{\mathbf{f}}}{\d t} &= (\mathbf{a}_\mathbf{0} + \mathbf{a}_\mathbf{1} \boldsymbol{\mu}_{\mathbf{f}}) + (\mathbf{R}_{\mathbf{f}}\mathbf{A}_\mathbf{1}^* ) (\mathbf{B}_1\mathbf{B}_1^*)^{-1} \left(\frac{\d\mathbf{X}}{\d t} - (\mathbf{A}_\mathbf{0} + \mathbf{A}_\mathbf{1}\boldsymbol{\mu}_{\mathbf{f}})\right),\label{CGNS_Stat_Mean}\\
  \frac{\d\mathbf{R}_{\mathbf{f}}}{\d t} &= \mathbf{a}_\mathbf{1} \mathbf{R}_{\mathbf{f}} + \mathbf{R}_{\mathbf{f}}\mathbf{a}_\mathbf{1}^* + \mathbf{b}_2\mathbf{b}_2^* - ( \mathbf{R}_{\mathbf{f}}\mathbf{A}_\mathbf{1}^*)(\mathbf{B}_1\mathbf{B}_1^*)^{-1}(\mathbf{A}_\mathbf{1}\mathbf{R}_{\mathbf{f}}),\label{CGNS_Stat_Cov}
\end{align}
\end{subequations}
with $\cdot^*$ being the complex conjugate transpose. The subscript `$\mathbf{f}$' in the conditional mean $\boldsymbol{\mu}_{\mathbf{f}}$ and conditional covariance $\mathbf{R}_{\mathbf{f}}$ is an abbreviation for `filter'. The explicit formulae in~\eqref{CGNS_PDF}-\eqref{CGNS_Stat} correspond to the optimal nonlinear filter solution of the state variable $\mathbf{Y}(t)$ given a realization of the observed time series $\mathbf{X}(s)$ for $s\in[0,t]$.
Thus,  $\boldsymbol\mu_{\mathbf{f}}$ and  $\mathbf{R}_{\mathbf{f}}$ in~\eqref{CGNS_Stat} are also known as the filter posterior mean and the filter posterior covariance. The classical Kalman-Bucy filter~\cite{kalman1961new} is the simplest special example of~\eqref{CGNS_Stat}.

The closed analytic formula~\eqref{CGNS_Stat} provides an efficient algorithm for the nonlinear data assimilation of the CGNS, which avoids using the ensemble or particle methods that may suffer from sampling errors. In Section~\ref{Sec:DA_Forecast}, the closed analytic data assimilation formula~\eqref{CGNS_Stat} will be used for an accurate state estimation of the initial value that facilitates effective ensemble forecast. It also allows the development of an efficient Gaussian mixture algorithm for calculating the non-Gaussian PDFs of the CGNS (see Section~\ref{Subsec: PDFs}), which overcomes the curse of dimensionality. Such non-Gaussian PDFs are crucial in analyzing the model sensitivity and predicting the system response, the details of which will be discussed in Section~\ref{Sec:Response}.

\subsubsection{Nonlinear smoother}\label{Subsec: smoother}
Filtering exploits the observational information up to the current time instant for an online state estimation. On the other hand, given the observational time series within an entire time interval, the state estimation can become more accurate. This is known as the smoother~\cite{sarkka2013bayesian}.

Given one realization of the observed variable $\mathbf{X}(t)$ for $t\in[0,T]$, the optimal smoother estimate $p(\mathbf{Y}(t)|\mathbf{X}(s), s\in[0,T])$ of the CGNS~\eqref{CGNS} is also Gaussian~\cite{chen2020learning},
\begin{equation}\label{Smoother}
  p(\mathbf{Y}(t)|\mathbf{X}(s), s\in[0,T])\sim\mathcal{N}(\boldsymbol\mu_\mathbf{s}(t),\mathbf{R}_\mathbf{s}(t)),
\end{equation}
where the conditional mean $\boldsymbol\mu_\mathbf{s}(t)$ and conditional covariance $\mathbf{R}_\mathbf{s}(t)$ of the smoother at time $t$ satisfy the following backward equations
\begin{subequations}\label{Smoother_Main}
\begin{align}
  \frac{\overleftarrow{\d \boldsymbol{\mu}_\mathbf{s}}}{\d t} &=  -\mathbf{a}_\mathbf{0} - \mathbf{a}_\mathbf{1}\boldsymbol{\mu}_\mathbf{s}  + (\mathbf{b}_2\mathbf{b}_2^*)\mathbf{R}_{\mathbf{f}}^{-1}(\boldsymbol\mu_{\mathbf{f}} - \boldsymbol{\mu}_\mathbf{s}),\label{Smoother_Main_mu}\\
  \frac{\overleftarrow{\d \mathbf{R}_\mathbf{s}}}{\d t} &= - (\mathbf{a}_\mathbf{1} + (\mathbf{b}_2\mathbf{b}_2^*) \mathbf{R}_{\mathbf{f}}^{-1})\mathbf{R}_\mathbf{s} - \mathbf{R}_\mathbf{s}(\mathbf{a}_\mathbf{1}^* + (\mathbf{b}_2\mathbf{b}_2^*)\mathbf{R}_{\mathbf{f}})  + \mathbf{b}_2\mathbf{b}_2^* ,\label{Smoother_Main_R}
\end{align}
\end{subequations}
with $\boldsymbol\mu_{\mathbf{f}}$ and $\mathbf{R}_{\mathbf{f}}$ given by~\eqref{CGNS_Stat}. Here, the subscript `$\mathbf{s}$' in the conditional mean $\boldsymbol{\mu}_{\mathbf{s}}$ and conditional covariance $\mathbf{R}_{\mathbf{s}}$ is an abbreviation for `smoother', which should not be confused by the time variable $s$ in $\mathbf{X}(s)$. The notation $\overleftarrow{\d \cdot}/\d t$ corresponds to the negative of the usual derivative, which means that the system \eqref{Smoother_Main} is solved backward over $[0,T]$ with $(\boldsymbol\mu_\mathbf{s}(T), \mathbf{R}_\mathbf{s}(T)) = (\boldsymbol\mu_{\mathbf{f}}(T), \mathbf{R}_{\mathbf{f}}(T))$ with the starting value of the nonlinear smoother $(\boldsymbol\mu_\mathbf{s}(T), \mathbf{R}_\mathbf{s}(T))$ being the same as the filter estimate $(\boldsymbol\mu_{\mathbf{f}}(T), \mathbf{R}_{\mathbf{f}}(T))$.

The nonlinear smoother plays an important role for an unbiased state estimation and postprocessing of the data. It is also able to quantify the uncertainty in the unobserved variables in the parameter estimation given only partial observations, which will be a topic to be studied in Section~\ref{Sec:Parameter_Estimation}. In addition, the nonlinear smoother is  the basis to the development of a nonlinear sampling formula, which will be shown below and is a necessary step in analyzing the model sensitivity and predicting the system response in Section~\ref{Sec:Response}.

\subsubsection{Nonlinear sampling formula}\label{Subsec: sampling}
Associated with the nonlinear smoother, a nonlinear conditional sampling formula can be derived. In addition to satisfying the point-wise optimal estimate~\eqref{Smoother_Main}, i.e., a distribution formed by conditional mean and conditional distribution at each time stamp, the conditional sampled trajectories further take into account the path-wise temporal correlation. These sampled trajectories in the CGNS framework can be regarded as the analogs of the ensemble members in the ensemble Kalman smoother~\cite{evensen2009data}, but the former can be obtained via a closed analytic formula.

Conditioned on one realization of the observed variable $\mathbf{X}(s)$ for $s\in[0,T]$, the optimal strategy of sampling  the trajectories associated with the unobserved variable $\mathbf{Y}$ satisfies the following explicit formula~\cite{chen2020efficient},
\begin{equation}\label{Sampling_Main}
  \frac{\overleftarrow{\d \mathbf{Y}}}{\d t} = \frac{\overleftarrow{\d \boldsymbol\mu_\mathbf{s}}}{\d t} - \big(\mathbf{a}_\mathbf{1} + (\mathbf{b}_2\mathbf{b}_2^*)\mathbf{R}_\mathbf{f}^{-1}\big)(\mathbf{Y} - \boldsymbol\mu_\mathbf{s}) + \mathbf{b}_2\dot{\mathbf{W}}_{\mathbf{Y}}(t),
\end{equation}
where $\dot{\mathbf{W}}_{\mathbf{Y}}(t)$ is a random noise that is independent from $\dot{\mathbf{W}}_{2}(t)$ in~\eqref{CGNS}. The conditional sampling formula is another necessary component in  analyzing the model sensitivity and predicting the system response in Section~\ref{Sec:Response}.

\subsection{Semi-analytic and statistically accurate formulae for solving the non-Gaussian PDFs via mixtures}\label{Subsec: PDFs}
The closed analytic formula~\eqref{CGNS_Stat} in calculating the conditional distribution $p(\mathbf{Y}(t)|\mathbf{X}(s),s\leq t)$ in~\eqref{CGNS_PDF} also provides an extremely useful way to compute the marginal distribution $p({\mathbf{Y}(t)})$. In fact, assuming there are $L$ trajectories of $\mathbf{X}(s\leq t)$, denoted by  $\mathbf{X}^{\text{obs}}_l(s\leq t)$ for $l=1,\ldots,L$, then in the limit with $L\to\infty$, the marginal distribution $p({\mathbf{Y}(t)})$ is given by
  \begin{equation}\label{Marginal_Y}
    p({\mathbf{Y}(t)}) =\lim_{L\to\infty}\frac{1}{L}\sum_{l=1}^L p\Big({\mathbf{Y}(t)}|{\mathbf{X}^{\text{obs}}_l(s\leq t)}\Big).
\end{equation}
While the above identity is in the asymptotic form with $L \rightarrow \infty$, it has been shown that the error bound in approximating $p({\mathbf{Y}(t)})$ with a finite $L$ does not depend on the dimension of $\mathbf{Y}$~\cite{chen2018rigorous}. In other words, fundamentally different from the traditional Monte Carlo simulations, the method in~\eqref{Marginal_Y} avoids the curse of dimensionality. If it is further assumed that the dimension of $\mathbf{X}$ is low, then the following efficient and statistically accurate approach can be utilized to compute the joint PDF at any transient phase $p(\mathbf{X}(t),\mathbf{Y}(t))$~\cite{chen2018efficient, chen2017beating},
\begin{equation}\label{Joint_X_Y}
    p(\mathbf{X}(t),\mathbf{Y}(t)) =\lim_{L\to\infty} \frac{1}{L}\sum_{l=1}^L \Big(K_\mathbf{H}(\mathbf{X}(t)-\mathbf{X}^{\text{obs}}_l(t))\, p(\mathbf{Y}(t)|\mathbf{X}^{\text{obs}}_l(s\leq t))\Big).
\end{equation}
Here the distribution of $\mathbf{X}$ is approximated by a mixture distribution that is solved via a kernel density estimation. One of the simplest choices, which is also the one used in this paper, is that each mixture component is a Gaussian function centered at $\mathbf{X}^{\text{obs}}_l(t)$. The same bandwidth $\mathbf{H}$ is utilized for different Gaussian mixture components and it is determined via the ``solve-the-equation'' method~\cite{botev2010kernel} that is an appropriate method for approximating non-Gaussian distributions.  In addition to the advantage of applying~\eqref{Joint_X_Y} in solving high-dimensional PDF (especially when the dimension of $\mathbf{Y}$ is large), the semi-analytic formula in~\eqref{Joint_X_Y} also allows a smoothed PDF, which reduces the sampling error compared with other approaches even for systems with moderate or low dimensions.

It is important to note that if the underlying system contains model error, then the PDF provided by~\eqref{Joint_X_Y} is very different from the one created by simply running the imperfect  model. This is because the model error in the marginal distribution $p({\mathbf{Y}(t)})$ is mitigated with the help from the observations $\mathbf{X}^{\text{obs}}_l(s\leq t)$. Such a unique feature plays a crucial role in improving the results for  computing the model response and sensitivity analysis, where the perfect model is seldom known in practice. The details will be illustrated in Section~\ref{Sec:Response}. As a final remark, only the equilibrium PDF is required in many applications, including the study of the model response. Therefore, if the system is ergodic, then only a single (sufficiently long) trajectory of $\mathbf{X}(0\leq t\leq T)$, denoted by $\mathbf{X}^{\text{obs}}(0\leq t\leq T)$, is needed in computing the equilibrium distribution $\peq(\mathbf{X},\mathbf{Y})$
\begin{equation}\label{Joint_X_Y_Equilibrium}
    \peq(\mathbf{X},\mathbf{Y}) =\lim_{J\to\infty} \frac{1}{J}\sum_{j=1}^J \Big(K_\mathbf{H}(\mathbf{X}-\mathbf{X}^{\text{obs}}(t_j))\, p(\mathbf{Y}|\mathbf{X}^{\text{obs}}(s\leq t_j))\Big),
\end{equation}
where $[t_1,\ldots,t_J]$ is a partition of the time interval $[T_0,T]$ with some burn-in time $T_0$.

\section{Strategies of Developing CGNS}\label{Sec:Model_Development}
The goals of developing approximate models for solving different problems in practice are often distinct to each other. For example, some applications require a skillful forecast model while others seek for a suitable model to reproduce certain statistics.
While there is not a universal criterion to build the ``optimal'' CGNS as an approximate model for all applications, a few potentially useful strategies are provided below for constructing CGNS, which will be applied in the following sections.

\subsection{Fast wave averaging}\label{Subsec:averaging}
Recall the general form of the complex systems with quadratic nonlinearity~\eqref{eq:abs_formu}. Writing it into the form of state variables $(\X,\Y)$ yields
\begin{equation}\label{Split_Form_General_Eqn}
\begin{split}
   \!\!\!      \frac{\mathrm{d}\X}{\mathrm{d}t} & \!=\! L_{11} \X  +  L_{12}\Y+B^{1}_{11} (\X,\!\X) +B^{1}_{12}(\X,\!\Y)+B^{1}_{22}(\Y,\!\Y) + \F_1(t)+\boldsymbol{\sigma}_1(\X,\!\Y, t)\dot{\mathbf{W}}_1(t),\\
   \!\!\! \frac{\mathrm{d}\Y}{\mathrm{d}t} & \!=\! L_{21}\X+L_{22}\Y+B^{2}_{11}(\X,\!\X) +B^{2}_{12}(\X,\!\Y)+B^{2}_{22}(\Y,\!\Y) +\F_2(t)+\boldsymbol{\sigma}_2(\X,\!\Y, t)\dot{\mathbf{W}}_2(t).
\end{split}
\end{equation}
In~\eqref{Split_Form_General_Eqn}, $L_{ij}$ are constant matrices while $B^{k}_{ij}$ are vector functions with the entries being quadratic functions of the state variables. In some applications, there exists a scale separation of the state variables, where $\X$ and $\Y$ represent slow and fast variables, respectively. In such a case, it is natural to apply a fast wave average such that the terms representing the self-interaction between the fast variables, i.e., $B^{1}_{22}\left(\Y,\Y\right)$ and $B^{2}_{22}\left(\Y,\Y\right)$, are approximated by stochastic damping and noise~\cite{majda2001mathematical, majda1999models}. By further approximating the diffusion coefficients $\boldsymbol{\sigma}_1\left(\X,\Y\right)$ and $\boldsymbol{\sigma}_2\left(\X,\Y\right)$ by functions of only $\X$, the resulting system becomes
\begin{equation}\label{Split_Form_CGNS}
\begin{split}
    \frac{\mathrm{d}\X}{\mathrm{d}t}&= \widetilde{L}_{11}\X+\widetilde{L}_{12}\Y+B^{1}_{11}\left(\X,\X\right)+B^{1}_{12}\left(\X,\Y\right)+\F_1(t)+\widetilde{\boldsymbol{\sigma}}_1\left(\X, t\right)\dot{\mathbf{W}}_1(t),\\
    \frac{\mathrm{d}\Y}{\mathrm{d}t}&= \widetilde{L}_{21}\X+\widetilde{L}_{22}\Y+B^{2}_{11}\left(\X,\X\right)+B^{2}_{12}\left(\X,\Y\right)+\F_2(t)+\widetilde{\boldsymbol{\sigma}}_2\left(\X, t\right)\dot{\mathbf{W}}_2(t),
\end{split}
\end{equation}
which belongs to the CGNS~\eqref{CGNS}. Note that, if the scale separation is not strong enough to apply the fast wave averaging, then $B^{1}_{22}\left(\Y,\Y\right)$ and $B^{2}_{22}\left(\Y,\Y\right)$ can be approximated by additional closure terms~\cite{san2018extreme, mou2021reduced} that include nonlinear functions of $\X$ and bilinear functions of $\X$ and $\Y$ to fit the CGNS framework.

\subsection{Stochastic parameterizations}\label{sec: SP}
The fast wave averaging or closure approximations are suitable approaches to build CGNS if the starting complex nonlinear system is completely known. However, in many practical applications, the information of the perfect model is not entirely available. Specifically, while the large-scale dynamics of nature is often accessible, the details of the small- or unresolved-scale features are not fully understood in many applications. In such a situation, suitable stochastic parameterizations can be adopted to approximate the processes of the unobserved variables $\Y$ such that the feedback from small/unresolved to large/resolved scales are well characterized and the parameterized system follows the CGNS structure.

One of the simplest strategies is to apply a linear stochastic model with Gaussian white noise to describe each component of the hidden processes of $\Y$ while the processes of $\X$ remain highly nonlinear. This parameterization strategy has been utilized in data assimilation and short-term statistical prediction~\cite{gershgorin2010test, gershgorin2010improving, branicki2013non}. The model with such a simple stochastic parameterization automatically fits the CGNS as there is no quadratic function of $\Y$ involved.
A more sophisticated stochastic parameterization is to incorporate the physics-constrained nonlinear regression model~\cite{majda2012physics, harlim2014ensemble} into the CGNS. This allows a more accurate way in characterizing the nonlinear dynamics of the small-scale features in $\Y$, influenced by the large-scale variables $\X$. In addition, the coupled system with physics constraints also prevents finite-time blow up of the solutions and facilitates a skillful medium- to long-range forecast.

\subsection{System augmentation} \label{sec:sys augmentation}
Another strategy to derive approximate CGNS is via a simple system augmentation technique detailed below. As has been discussed in Section~\ref{Subsec:averaging} that the most significant difference between the general nonlinear system \eqref{eq:abs_formu} and the CGNS \eqref{CGNS} is the quadratic nonlinear self-interactions of $\Y$, namely the terms $B^{1}_{22}\left(\Y,\Y\right)$ and $B^{2}_{22}\left(\Y,\Y\right)$ appearing in \eqref{Split_Form_General_Eqn}. After suitable regrouping of the terms in \eqref{Split_Form_General_Eqn} and assuming for simplicity that the diffusion coefficients depend only on $X$, we can rewrite this system \eqref{Split_Form_General_Eqn} into the following form,
\begin{subequations}\label{CGNS_2}
\begin{align}
  \frac{\d\mathbf{X}}{\d t} &= \Big[\mathbf{A}_\mathbf{0}(\mathbf{X},t) + \mathbf{A}_\mathbf{1}(\mathbf{X},t) \mathbf{Y} + B^{1}_{22}(\mathbf{Y},\mathbf{Y})\Big]  + \boldsymbol{\sigma}_1(\mathbf{X},t)\dot{\mathbf{W}}_\mathbf{1}(t),\label{CGNS_X_2}\\
  \frac{\d\mathbf{Y}}{\d t} &= \Big[\mathbf{a}_\mathbf{0}(\mathbf{X},t) + \mathbf{a}_\mathbf{1}(\mathbf{X},t) \mathbf{Y} + B^{2}_{22}(\mathbf{Y},\mathbf{Y})\Big]   + \boldsymbol{\sigma}_2(\mathbf{X},t)\dot{\mathbf{W}}_\mathbf{2}(t),\label{CGNS_Y_2}
\end{align}
\end{subequations}
which differs from the CGNS \eqref{CGNS} by the two quadratic terms  $B^{1}_{22}(\mathbf{Y},\mathbf{Y})$ and $B^{2}_{22}(\mathbf{Y},\mathbf{Y})$. Instead of approximating these two terms directly via fast wave averaging as proposed in Section~\ref{Subsec:averaging}, we consider here the situation that the scale separation between $\mathbf{X}$ and $\mathbf{Y}$ is not pronounced.

We will still handle $B^{2}_{22}(\mathbf{Y},\mathbf{Y})$ in \eqref{CGNS_Y_2} via suitable stochastic parameterization in terms of $\mathbf{X}$, leading to an approximate equation for $\mathbf{Y}$ of the form
\begin{equation}
\frac{\d\mathbf{Y}}{\d t} = \Big[\widetilde{\mathbf{a}}_\mathbf{0}(\mathbf{X},t) + \widetilde{\mathbf{a}}_\mathbf{1}(\mathbf{X},t) \mathbf{Y}\Big]   + \widetilde{\boldsymbol{\sigma}}_2(\mathbf{X},t)\dot{\mathbf{W}}_\mathbf{2}(t).\label{CGNS_Y_2b}
\end{equation}
The term $B^{1}_{22}(\mathbf{Y},\mathbf{Y})$ is then dealt with through a system augmentation strategy as explained below.
Denote by $\mathbf{Z} = ((Y_1)^2, Y_1Y_2, Y_1Y_3, \ldots)^\top$ the $N_2(N_2 + 1)/2$-dimensional vector whose components consist of all possible quadratic monomials involving the components of the $N_2$-dimensional small-scale variable $\mathbf{Y}$. Using It\^o's formula \cite{gardiner2004handbook} and \eqref{CGNS_Y_2b}, we can obtain the corresponding equation for $\mathbf{Z}$, which takes the following form
\begin{equation}\label{Variable_Z}
\begin{gathered}
  \frac{\d\mathbf{Z}}{\d t} = \Big[\mathbf{c}_\mathbf{0}(\mathbf{X},t) + \mathbf{c}_\mathbf{1}(\mathbf{X},t) \mathbf{Y} + \mathbf{c}_\mathbf{2}(\mathbf{X},t) \mathbf{Z} \Big] + \boldsymbol{\sigma}_3(\mathbf{X},\mathbf{Y},t)\dot{\mathbf{W}}_\mathbf{2}(t).
\end{gathered}
\end{equation}
Note that the quadratic term $B^{1}_{22}(\mathbf{Y},\mathbf{Y})$ in \eqref{CGNS_X_2} can be rewritten as a linear function of $\mathbf{Z}$, denoted by $L \mathbf{Z}$, thanks to the very choice of $\mathbf{Z}$. Thus, if we further approximate $\mathbf{Y}$ in the diffusion coefficient $\boldsymbol{\sigma}_3(\mathbf{X},\mathbf{Y},t)$ of \eqref{Variable_Z} by, e.g., its global mean $\overline{\mathbf{Y}}$, the augmented system for $(\mathbf{X}, \mathbf{Y}, \mathbf{Z})$ consisting of \eqref{CGNS_X_2}, \eqref{CGNS_Y_2b}, and \eqref{Variable_Z} fits into the form of CGNS given by \eqref{CGNS} with $(\mathbf{Y}, \mathbf{Z})$ here playing the role of $\mathbf{Y}$ in \eqref{CGNS}, after replacing $B^{1}_{22}(\mathbf{Y},\mathbf{Y})$ in \eqref{CGNS_X_2} by $L \mathbf{Z}$ and approximating $\boldsymbol{\sigma}_3(\mathbf{X},\mathbf{Y},t)$ in \eqref{Variable_Z} by $\boldsymbol{\sigma}_3(\mathbf{X},\overline{\mathbf{Y}},t)$.

With the CGNS structure available through this augmented system, the conditional distribution such as $p(\mathbf{Y}(t)|\mathbf{X}(s\leq t))$ can now be approximated by first computing the conditional distribution $p(\mathbf{Y}(t),\mathbf{Z}(t)|\mathbf{X}(s\leq t))$ for the augmented system using formulae provided in Section~\ref{Subsec:Pathwise} and then marginalizing it over $\mathbf{Z}$. The essence of this simple approach is thus to increase the dimension of the system in exchange of analytic formulae of the conditional distributions, which would otherwise be computationally expensive to compute even for models with moderate dimensions.

When the dimension, $N_2$, of $\mathbf{Y}$ is very high, one may prefer to identify only a subset of $\mathbf{Y}$ to apply the system augmentation technique in order not to inflate too much the number of variables, since the auxiliary variable $\mathbf{Z}$ has dimension $N_2(N_2 + 1)/2$. Such an extension goes beyond the scope of the current article, and will be addressed in a separate communication. Note also that if $B^{2}_{22}(\mathbf{Y},\mathbf{Y})\equiv0$, then Eq.~\eqref{CGNS_Y_2b} is reduced to \eqref{CGNS_Y_2}. In this case, there is no approximation involved in the drift part of \eqref{Variable_Z}.

We demonstrate now the efficiency of this approach in the context of data assimilation and ensemble forecast using a low-dimensional truncation of a stochastic Burgers-type equation.

\section{Data Assimilation and Ensemble Forecast}\label{Sec:DA_Forecast}

Data assimilation concerns the problem of estimating the state variables of a given, usually nonlinear and possibly stochastic, dynamical system when observations of certain related output variables are available \cite{evensen2009data, kalnay2003atmospheric, law2015data, majda2012filtering}. One major challenge in data assimilation is the strong nonlinearity and the associated non-Gaussian statistics in the underlying dynamics, in which a direct application of the particle methods may be inaccurate especially in the high dimensional situations.
The development of cheap and effective approximate models that capture the main characteristics of the underlying dynamics is thus an important topic in state estimation and data assimilation. Since the data assimilation solution corresponds to the initialization of the subsequent forecast, an efficient and accurate data assimilation scheme is also essential to advancing the forecast skill. Note that there is usually a stronger demand in developing suitable approximate models for data assimilation than the subsequent short- or medium-range forecast since the former often involves many numerical or sampling issues in the presence of strong nonlinearity and non-Gaussianity.

In this section, a particular type of approximate models for this purpose obtained by the method of system augmentations presented in Section~\ref{sec:sys augmentation} is studied. The resulting approximate model has the form of a CGNS. Thus, the associated data assimilation solutions can be calculated using the closed analytic formula \eqref{CGNS_Stat} as was discussed in Section~\ref{Subsec:Pathwise}. To simplify the presentation, the idea is illustrated using a low-dimensional stochastic differential equation (SDE) with energy-conserving quadratic terms. The data assimilation results from the CGNS, which is an approximate model, will be compared with that by applying a classical ensemble data assimilation method directly to the perfect model. The goal is to illustrate the efficiency and accuracy of the data assimilation scheme using the CGNS, especially in avoiding the potential sampling and other numerical issues that appear in the ensemble-based approaches.

\subsection{A truncated stochastic quadratic system and its CGNS approximation through system augmentation} \label{sec: TBH model}
The model considered here is the following three-dimensional SDEs with energy-conserving quadratic nonlinear terms and subject to additive white noise forcing
\begin{subequations} \label{Eq_xyz}
\begin{align}
& \frac{\d x}{\d t} =  \beta_x x + \alpha x y + \alpha  yz   + \sigma_x \dot{W}_x, \label{Eq1_xyz}\\
& \frac{\d y}{\d t} =  \beta_y y - \alpha x^2 + 2\alpha xz  + \sigma_y \dot{W}_y, \label{Eq2_xyz} \\
& \frac{\d z}{\d t} = \beta_z z - 3 \alpha x y + \sigma_z \dot{W}_z. \label{Eq3_xyz}
\end{align}
\end{subequations}
Here, the coefficients for the linear terms are chosen such that $\beta_x$ is positive to introduce linear instability into the system, while $\beta_y$ and $\beta_z$ are negative, representing linear damping effects. The coefficient $\alpha>0$ controls the strength of the nonlinearity; and the noise strength coefficients $\sigma_x$, $\sigma_y$, and $\sigma_z$ are positive constants. This system can for instance be obtained as a Fourier-Galerkin projection of the stochastic Burgers-Sivashinsky equation
\begin{equation*}
\frac{\partial u}{\partial t}  = \big( \nu \partial_{xx} u  + \lambda u  - u  \partial_x u\big) + \dot{W}(t,x)
\end{equation*}
posed on a bounded interval $x \in (0,L)$ subject to homogeneous Dirichlet boundary conditions. In this context, $\beta_x$, $\beta_y$, and $\beta_z$ are simply the three largest eigenvalues of the linear operator and $\alpha$ is linked to the domain size $L$ via $\alpha = \pi/(\sqrt{2} L^{3/2})$; see e.g. \cite[Chapter 6]{CLW15_vol2}.

In the following, the largest scale variable $x$ is treated as the observed variable while there is no direct observations for the state variables $(y,z)$. Under this splitting of the state variables, system \eqref{Eq_xyz} does not have the conditional Gaussian structure due to the quadratic nonlinear term $\alpha yz$ between the unobserved variables that appears in \eqref{Eq1_xyz}.
Following the idea presented in Section~\ref{sec:sys augmentation}, in order to obtain a CGNS to approximate the system \eqref{Eq_xyz}, three auxiliary variables are introduced for the possible quadratic interactions between the two unobserved variables:
\begin{equation}
p = y^2, \quad q = yz, \quad r = z^2.
\end{equation}
Using \eqref{Eq_xyz} and apply It\^o's formula yields
\begin{equation} \label{Eq_pqr}
\begin{split}
& \frac{\d p}{\d t} = (\sigma_y)^2 +  2 \big( \beta_y p - \alpha x^2y + 2\alpha xq \big)   + 2\sigma_y y \dot{W}_y, \\
& \frac{\d q}{\d t}  = (\beta_y + \beta_z)q - \alpha x^2 z - 3 \alpha xp + 2\alpha xr    + \sigma_y z \dot{W}_y +   \sigma_y y \dot{W}_z, \\
& \frac{\d r}{\d t} = (\sigma_z)^2 + 2 \big( \beta_z r - 3 \alpha x q\big)  + 2 \sigma_z z \dot{W}_z.
\end{split}
\end{equation}

Assume that the global mean values of the unobserved variables $y$ and $z$ are accessible (from a period of training data). Then, in combination with \eqref{Eq_xyz}, and after replacing $y$ and $z$ in the state-dependent noise terms of \eqref{Eq_pqr} by their respective global mean, the following augmented system is arrived at:
\begin{equation} \label{Eq_augmented}
\begin{split}
& \frac{\d x}{\d t} =  \beta_x x + \alpha x y + \alpha  q  + \sigma_x \dot{W}_x, \\
& \frac{\d y}{\d t} =  \beta_y y - \alpha x^2 + 2\alpha xz + \sigma_y \dot{W}_y, \\
& \frac{\d z}{\d t} =  \beta_z z - 3 \alpha x y  + \sigma_z \dot{W}_z,\\
& \frac{\d p}{\d t} = (\sigma_y)^2 +  2 \big( \beta_y p - \alpha x^2y + 2\alpha xq \big)   + 2\sigma_y \overline{y} \dot{W}_y, \\
& \frac{\d q}{\d t}  = (\beta_y + \beta_z)q - \alpha x^2 z - 3 \alpha xp + 2\alpha xr    + \sigma_y \overline{z} \dot{W}_y +   \sigma_z \overline{y} \dot{W}_z, \\
& \frac{\d r}{\d t} = (\sigma_z)^2 + 2 \big( \beta_z r - 3 \alpha x q\big)  + 2 \sigma_z \overline{z} \dot{W}_z,
\end{split}
\end{equation}
where $yz$ in \eqref{Eq1_xyz} becomes the state variable $q$. This augmented system \eqref{Eq_augmented} fits into the CGNS form of \eqref{CGNS} with now the unobserved variables taken to be  $\Y = (y, z, p, q, r)^\top$.

Although the dimension of the approximate system is  increased compared with the original system,  closed analytic equations are now accessible for the evolution of the corresponding conditional statistics for the data assimilation solutions (see equation \eqref{CGNS_Stat} in Section~\ref{Subsec:Pathwise}).
As will be shown below, the approximate system \eqref{Eq_augmented} can provide a significantly more accurate estimation of $(y,z)$ compared with another conditional Gaussian approximation obtained by simply removing the term $\alpha yz$ in \eqref{Eq1_xyz}, called the bare truncation (BT) system below. The skill of the proposed method is comparable and sometimes even more accurate than the ensemble Kalman-Bucy filter (EnKBF) \cite{bergemann2012ensemble}, while being more efficient thanks to the availability of analytic formulae.

\subsection{Dynamical regimes and numerical setup} \label{sec: TBH numerical_setup}
In the following, we consider two dynamical regimes:
\begin{equation} \label{Eq_TBH_regimes}
\begin{split}
\text{Regime I:} & \quad \sigma_x = 1, \sigma _y = 1, \sigma_z = 2, \beta_x = 0.1, \beta_y = -0.5, \beta_z = -1, \alpha = \pi/\sqrt{2}, \\
\text{Regime II:} & \quad \sigma_x = 0.1, \sigma _y = 1, \sigma_z = 2, \beta_x = 0.1, \beta_y = -0.5, \beta_z = -1, \alpha = \pi/\sqrt{2}.
\end{split}
\end{equation}
In particular, we have a relatively strong nonlinear effects with $\alpha \approx 2.2$, and a relatively small spectral gap between the observed and the hidden variables with $\beta_x - \beta_y = 0.6$. Both of the two hidden variables are subject to strong noise perturbations. The two regimes differ only in the value of the noise strength $\sigma_x$ in the $x$-equation.

The same numerical setup is adopted for both parameter regimes. The true signal is obtained by integrating the original SDE system \eqref{Eq_xyz} for an arbitrarily fixed noise path using the Euler-Maruyama scheme with a uniform time step size $\delta t = 5\times 10^{-4}$ and initialized at $(x,y,z) = (0,0,0)$.

For the state estimation of the unobserved variables $(y,z)$, we compare three methods:

\begin{itemize}
\item[] {\bf Method 1}: Apply the nonlinear filtering formulae \eqref{CGNS_Stat} for the general CGNS \eqref{CGNS} to the augmented system \eqref{Eq_augmented}, with $\X = x$ and $\Y = (y, z, p, q, r)^\top$. This method will be referred as the CG method below.

\item[] {\bf Method 2}: Apply the nonlinear filtering formulae \eqref{CGNS_Stat} to a bare truncation of \eqref{Eq_xyz} in which we simply remove the term $\alpha y z$ in \eqref{Eq1_xyz} to obtain a CGNS, with $\X = x$ and $\Y = (y, z)^\top$. This method will be referred as the BT method below.

\item[] {\bf Method 3}: Apply the ensemble Kalman-Bucy filtering (EnKBF) method to \eqref{Eq_xyz}.  See \eqref{Eq_EnKBF} below for its formulation.
\end{itemize}

The data assimilation for each of the above methods is performed over the time window $[0, 400]$ with the same time step size $\delta t$ as the true signal. For the CG method, the global mean values $\overline{y}$ and $\overline{z}$ in \eqref{Eq_augmented} are taken to be the mean values of the corresponding true signal over the interval $[0,200]$. For both CG and BT, the initial values of the conditional mean and conditional covariance are taken to be zero. For EnKBF, the size of ensemble is taken to be $N=100$ and the unobserved variables are initialized at $(y,z) = (0,0)$. For the sake of clarity, we provide below some details about the EnKBF applied to \eqref{Eq_xyz}. We introduce the following notations for the drift part of the system \eqref{Eq_xyz}:
\begin{equation}
\begin{split}
& g(x,y,z) =  \beta_x x + \alpha x y + \alpha  yz, \\
& f_1(x,y,z) = \beta_y y - \alpha x^2 + 2\alpha xz, \\
& f_2(x,y,z) =  \beta_z z - 3 \alpha x y.
\end{split}
\end{equation}
Denote by $\y = (y_1, y_2, \ldots, y_N)^\top$ and $\z = (z_1, z_2, \ldots, z_N)^\top$ the collection of all the $N$ ensemble members. We define also
\begin{equation}\label{EnKBF_Variance}
\begin{split}
& \mathcal{N}_{1}(x_{\text{obs}}(t),\y,\z) = \frac{1}{\sigma_x^2 (N-1)} \sum_{j=1}^N (y_j - \overline{\y}(t))(g(x_{\text{obs}}(t), y_j, z_j) - \overline{g}(x_{\text{obs}}(t), \y, \z)), \\
& \mathcal{N}_{2}(x_{\text{obs}}(t),\y,\z) = \frac{1}{\sigma_x^2 (N-1)} \sum_{j=1}^N (z_j - \overline{\z}(t))(g(x_{\text{obs}}(t), y_j, z_j) - \overline{g}(x_{\text{obs}}(t), \y, \z)), \\
\end{split}
\end{equation}
where
\begin{equation}
\overline{\y}(t) = \frac{1}{N}\sum_{\ell=1}^N  y_\ell(t), \quad \overline{\z}(t) = \frac{1}{N}\sum_{\ell=1}^N  z_\ell(t), \quad  \overline{g}(x_{\text{obs}}(t), \y, \z) = \frac{1}{N}\sum_{\ell=1}^N g(x_{\text{obs}}(t), y_{\ell}, z_{\ell}).
\end{equation}
Then, each ensemble member $(y_i, z_i)$, $i=1,2,\ldots,N$, of the EnKBF is computed using
\begin{equation} \label{Eq_EnKBF}
\begin{split}
\frac{\d y_i}{\d t} & = f_1(x_{\text{obs}}(t),y_i,z_i) + \sigma_y \dot{W}_{y,i}\\
&\qquad\qquad\qquad - \mathcal{N}_1(x_{\text{obs}}(t),\y,\z) \big[g(x_{\text{obs}}(t),y_i,z_i) -  \dot{x}_{\text{obs}}(t) + \sigma_x \dot{W}_{x,i} \big], \\
\frac{\d z_i}{\d t} & = f_2(x_{\text{obs}}(t),y_i,z_i) + \sigma_z \dot{W}_{z,i}\\
&\qquad\qquad\qquad - \mathcal{N}_2(x_{\text{obs}}(t),\y,\z) \big[g(x_{\text{obs}}(t),y_i,z_i) -  \dot{x}_{\text{obs}}(t) + \sigma_x \dot{W}_{x,i} \big],
\end{split}
\end{equation}
where $W_{x,i}$, $W_{y,i}$, and $W_{z,i}$, $i=1,2,\ldots,N$, are all mutually independent one-dimensional Brownian motions, and $x_{\text{obs}}$ is the observed signal of $x$.

For Regime II, we will also compare the ensemble forecast skills. The forecast is performed over the time window $[200,400]$, which is chosen to avoid overlap with the training window $[0,200]$ from which the global mean values of $y$ and $z$ appearing in \eqref{Eq_augmented} are computed. The forecast model is taken to be the true SDE system \eqref{Eq_xyz}, and the initial conditions (IC) of $(y,z)$ are drawn from multivariate Gaussian distributions with mean and covariance estimated respectively from BT, CG and EnKBF described above. For $x$, its initial value is taken to be that of the true signal at the corresponding time instant. We will also compute the results when the forecast is initialized with the true signal for all the three variables, which serves as the reference of the theoretic forecast/predictability limit and will be referred as the case with the perfect IC. The time locations at which to issue the forecasts are equally spaced over the chosen time interval, with a gap of 0.01 between two adjacent forecasts, leading thus to a total of $2\times 10^4$ forecast locations. Each forecast is computed up to a lead time of 1 time unit, and a total of 40 ensemble members are generated at each forecast location. This procedure is repeated for each of the methods used to construct the IC.

\subsection{Numerical results} \label{sec: TBH results}

We present now the results obtained based on the numerical procedure described above. For the two regimes given by \eqref{Eq_TBH_regimes},  due to the larger noise strength parameter $\sigma_x$ used in Regime I for the observed variable, the corresponding DA exercise is less challenging and will be presented first.

\begin{figure}[tbh!]
  \centering
  \includegraphics[width=\textwidth]{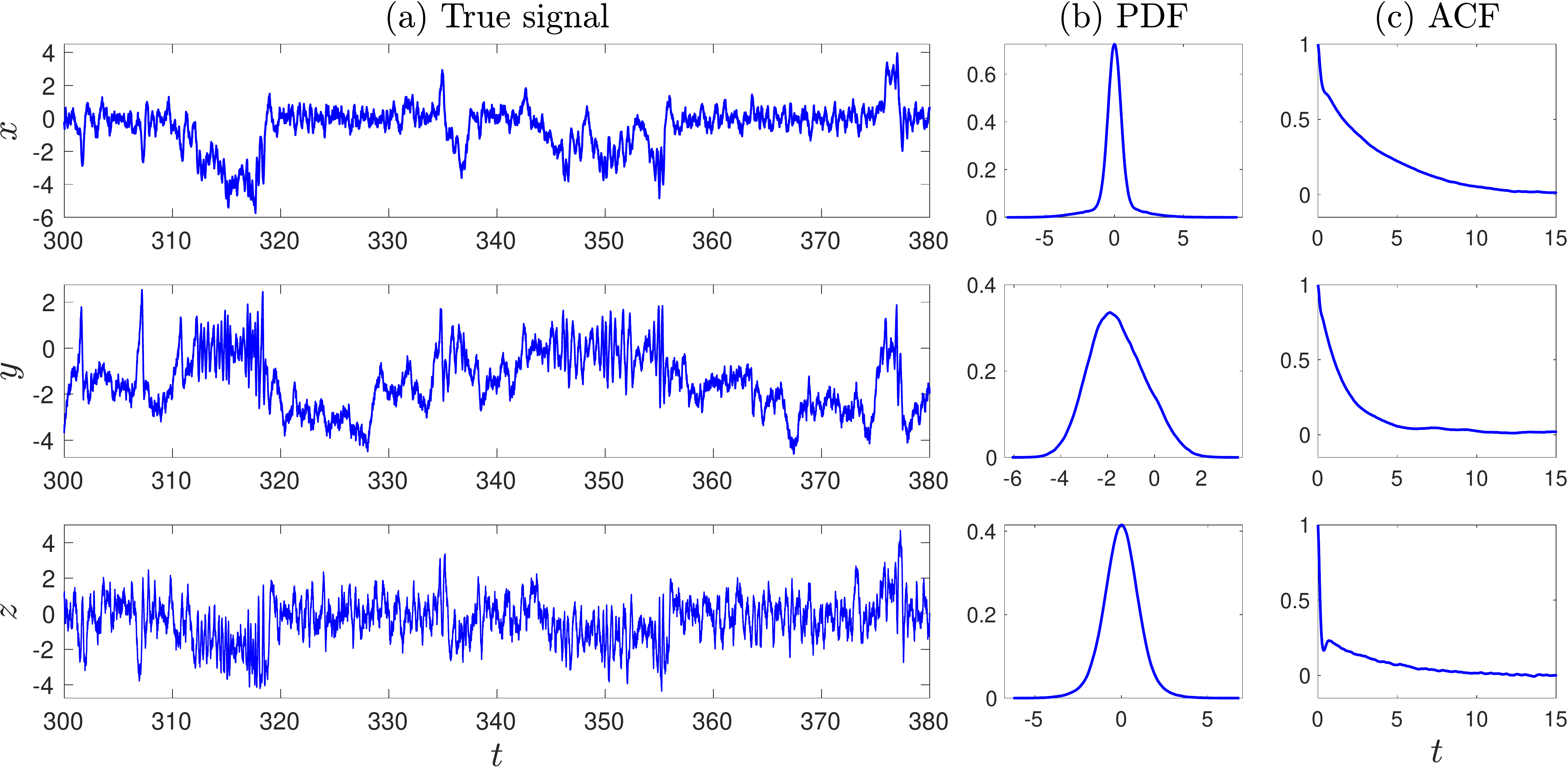}
  \caption{Solution of \eqref{Eq_xyz} in Regime I given by \eqref{Eq_TBH_regimes} for one arbitrarily fixed realization of the noise, and the corresponding probability density functions (PDFs) and autocorrelation functions (ACFs). This solution trajectory is taken to be the true signal. See Section~\ref{sec: TBH numerical_setup} for details about the numerical setup. The ACFs and PDFs are estimated based on the solution trajectory over the time window $[0,10^{4}]$, corresponding to $2\times 10^7$ data points for the time step used.}\label{fig: TBH_truth_case1}
\end{figure}

\medskip
\noindent{\bf Results for Regime I}. As is shown in Figure~\ref{fig: TBH_truth_case1} for Regime I, the dynamics of $x$ exhibits intermittent behavior with relatively quiescent  episodes punctuated by large excursion events. Due to the relatively small spectral gap, the dynamics of $y$ also exhibits highly nonlinear oscillations sustained by noise. In contrast, the dynamics of $z$ is mainly a damped oscillation sustained by noise due to the relatively strong linear stabilizing effects, and it is the variable that decays the fastest.

\begin{figure}[tbh!]
  \centering
  \includegraphics[width=0.8\textwidth]{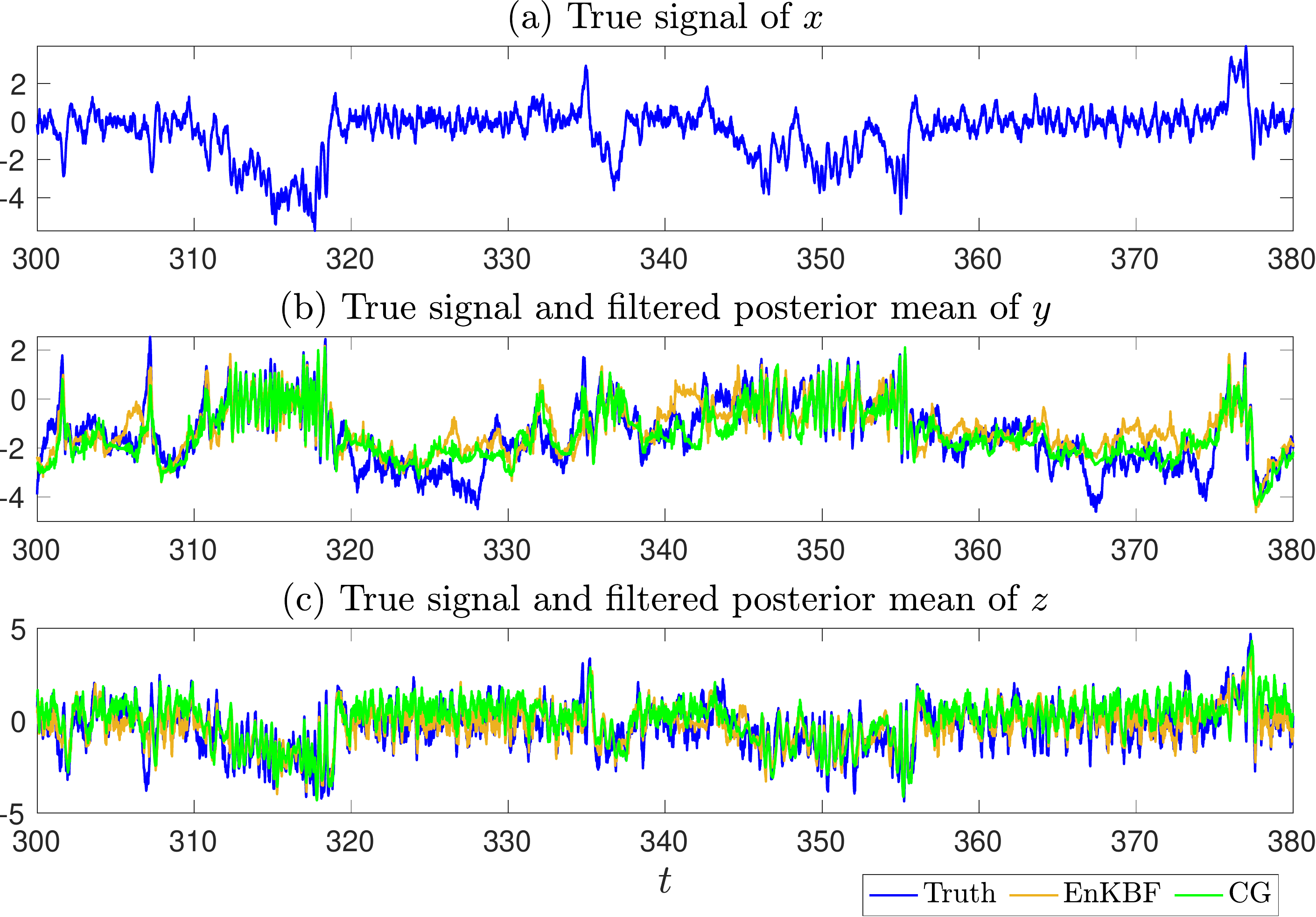}
  \caption{Panel (b): the filtered posterior mean of $y$ for Regime I obtained from CG (green) and EnKBF (orange); the corresponding true signal previously shown in Figure~\ref{fig: TBH_truth_case1} is plotted in blue. Panel (c): analogue of Panel (b) for $z$.}\label{fig: TBH_DA_case1}
\end{figure}

\begin{figure}[bth!]
  \centering
  \includegraphics[width=0.8\textwidth]{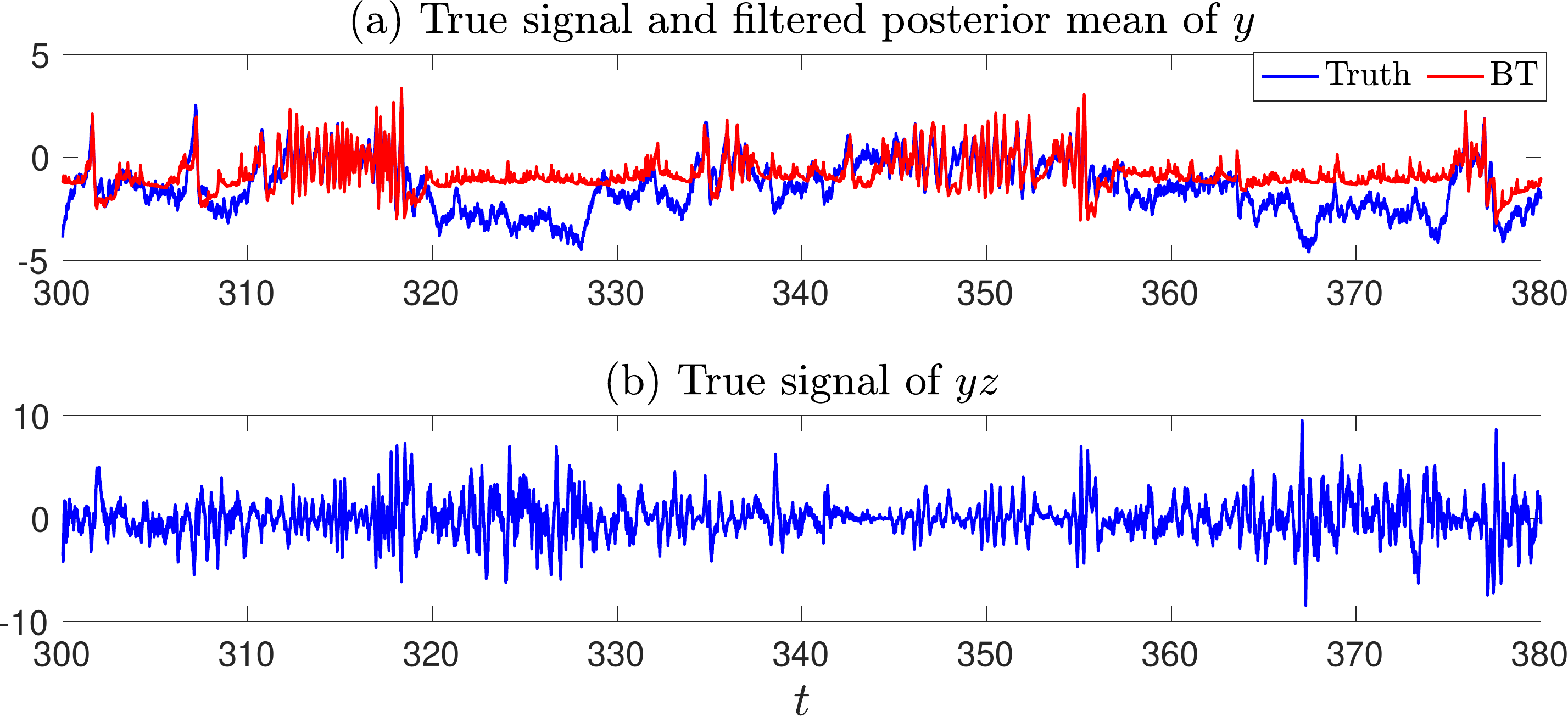}
  \caption{Panel (a): The filtered posterior mean of $y$ obtained from BT (red) and the true signal of $y$ (blue). Panel (b): The true signal of $yz$. The deterioration of the estimated $y$ using BT occurs over time windows when the magnitude of $yz$ gets large.}\label{fig: TBH_DA_BT_case1}
\end{figure}

The posterior mean states of $(y,z)$ for Regime I obtained by CG and EnKBF are fairly close to each other as is shown in Figure~\ref{fig: TBH_DA_case1}. The conditional covariance matrices of $(y,z)$ estimated by these two methods are also close to each other for this regime (not shown). In contrast, for BT, there are prolonged time windows over which the posterior mean of $y$ deviates significantly from the true signal as is shown in Panel (a) of Figure~\ref{fig: TBH_DA_BT_case1}. An inspection of the time series of $yz$ shown in Panel (b) of Figure~\ref{fig: TBH_DA_BT_case1} reveals that such deviation typically occurs when the value of $yz$ is large, which is expected, since the omission of the term $\alpha yz$ in \eqref{Eq1_xyz} is the only difference between the BT system and the full system. The posterior mean of $z$ obtained by BT is similar to those obtained by CG and EnKBF shown in Panel (c) of Figure~\ref{fig: TBH_DA_case1}, which is thus not presented.

\medskip
\noindent{\bf Results for Regime II}. The dynamics of the true system \eqref{Eq_xyz} in Regime II exhibits similar features as in Regime I shown in Figures~\ref{fig: TBH_truth_case1}, although the amplitude of each variable is slightly reduced due to the smaller noise intensity $\sigma_x$ used for this regime. The shape of both the PDFs and ACFs of all the three variables are similar to those shown in Figures~\ref{fig: TBH_truth_case1}, except that the decorrelation time of $x$ becoming comparable with that of $y$ in this regime. The analogue of Figures~\ref{fig: TBH_truth_case1} for Regime II is thus omitted.

\begin{figure}[tbh!]
  \centering
  \includegraphics[width=\textwidth]{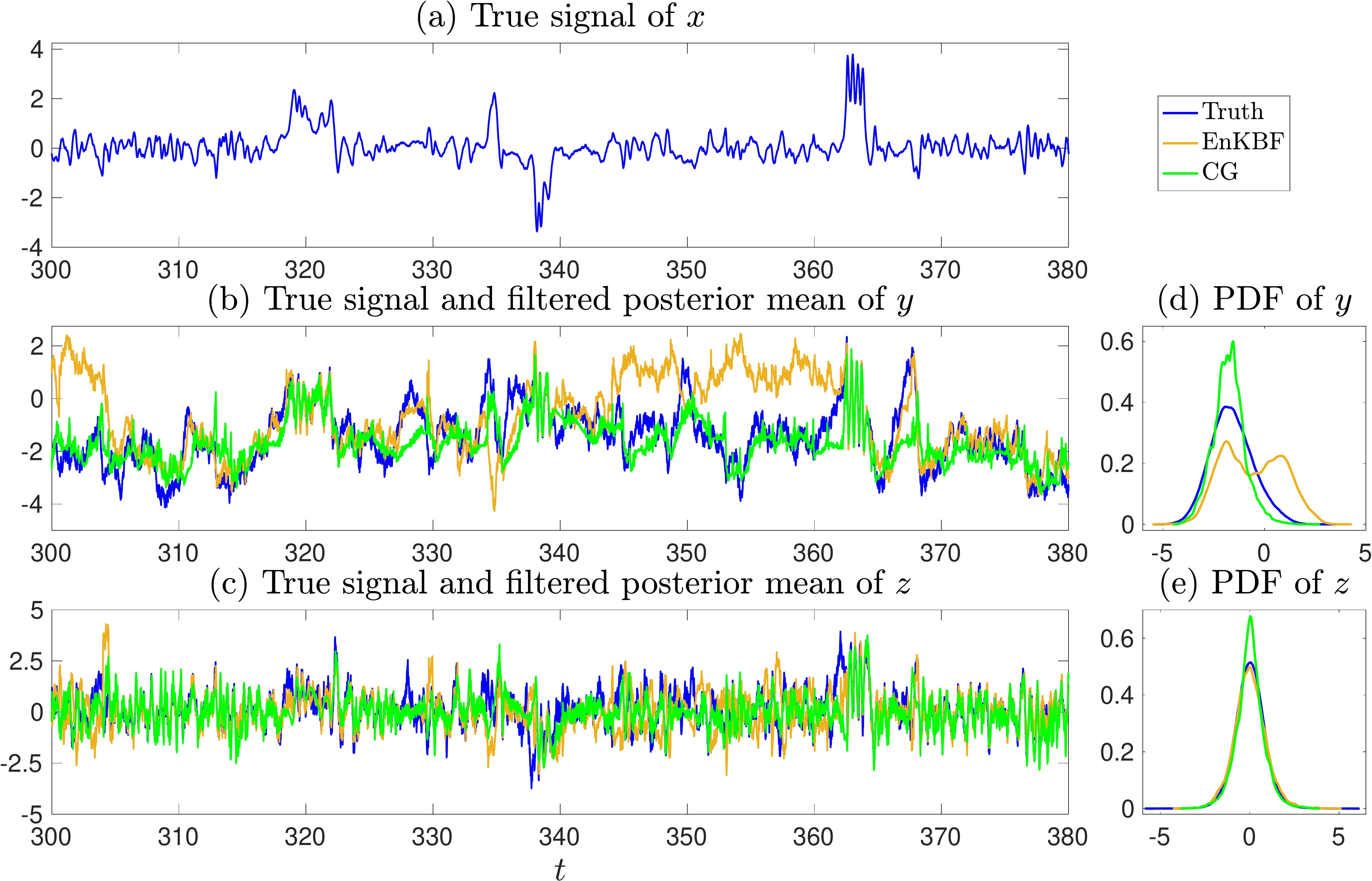}
  \caption{Panel (b): the filtered posterior mean of $y$ for Regime II obtained from CG (green) and EnKBF (orange) with the corresponding true signal shown in blue. Panel (c): the analogue of Panel (b) for $z$. Panel (d): PDFs of the filtered posterior mean of $y$ from CG (green) and EnKBF (orange) compared with that of the true signal. Panel (e): the analogue of Panel (d) for $z$.}\label{fig: TBH_DA_case2}
\end{figure}

For this regime, CG and EnKBF still provide comparable posterior mean state of $z$ as shown in Panel (c) of Figure~\ref{fig: TBH_DA_case2}, although the PDF of the posterior mean obtained by EnKBF approximates slightly better the PDF of the true signal of $z$. However, CG is significantly more skillful in estimating the conditional mean of $y$ (Figure~\ref{fig: TBH_DA_case2}, Panels (b) and (d)).

The deterioration of the skill from EnKBF in this regime is associated with a false bimodal behavior appearing in the posterior mean of the $y$ variable (see Panel (d) of Figure~\ref{fig: TBH_DA_case2}), whereas the true signal is unimodal, skewed towards negative values. It has also been checked that the bimodality is always there for EnKBF by further increasing the total number of ensemble members $N$ or decreasing the numerical integration time step $\delta t$. Such a pathological behavior is associated with the filter divergence \cite{gottwald2013mechanism, kelly2015concrete}, which often occurs for the ensemble-based filters when the  noise in the observational process is small and the observational process is highly nonlinear. In fact, the small  noise in the observational process $x$ makes the filter trusts more towards the information provided by the observations. However, the strong nonlinear and non-Gaussian features of $x$ make it very difficult to accurately recover the states of $y$ and $z$ by inferring mainly from the $x$ process.
In contrast, CG tracks well the modulations of the true signal, leading to a much better reproduction of the PDF of the true signal of $y$; see again Panel (d) of Figure~\ref{fig: TBH_DA_case2}. It is worth pointing out that the original system \eqref{Eq_xyz} can also exhibit bimodal dynamics in a broad range of dynamical regimes, even though bimodality is not observed in the true signals of $(x,y,z)$ for neither of the two regimes considered here. This bimodality that can occur in the dynamics of \eqref{Eq_xyz} is induced by the additive noise that drives the system to switch from the two locally stable steady states of the corresponding deterministic system produced from a supercritical pitchfork bifurcation, although when occurs, the bimodality is mainly visible in the $x$ variable.

Regarding BT, compared with the corresponding result shown in Figure~\ref{fig: TBH_DA_BT_case1} for Regime I, its performance here is even worse and is thus not shown. In particular, the posterior mean state of $y$ not only deviates significantly from the true signal but also has spurious fast oscillations presenting throughout the whole time window. Such fast oscillations also appear in the filtered posterior mean of $z$, although to a lesser extent.

For this regime, we also compared the skills of ensemble forecast with the forecast model simply taken to be the true SDE system \eqref{Eq_xyz}, and the initial conditions (IC) of $(y,z)$ drawn from multivariate Gaussian distributions in which the mean and covariance are estimated respectively from BT, CG and EnKBF described above; see Section~\ref{sec: TBH numerical_setup} for further details. In addition, we compute the results when the forecast is initialized with the true signal for all the three variables, which serves as the reference of the theoretic forecast limit.

\begin{figure}[tbh!]
  \centering
  \includegraphics[width=0.8\textwidth]{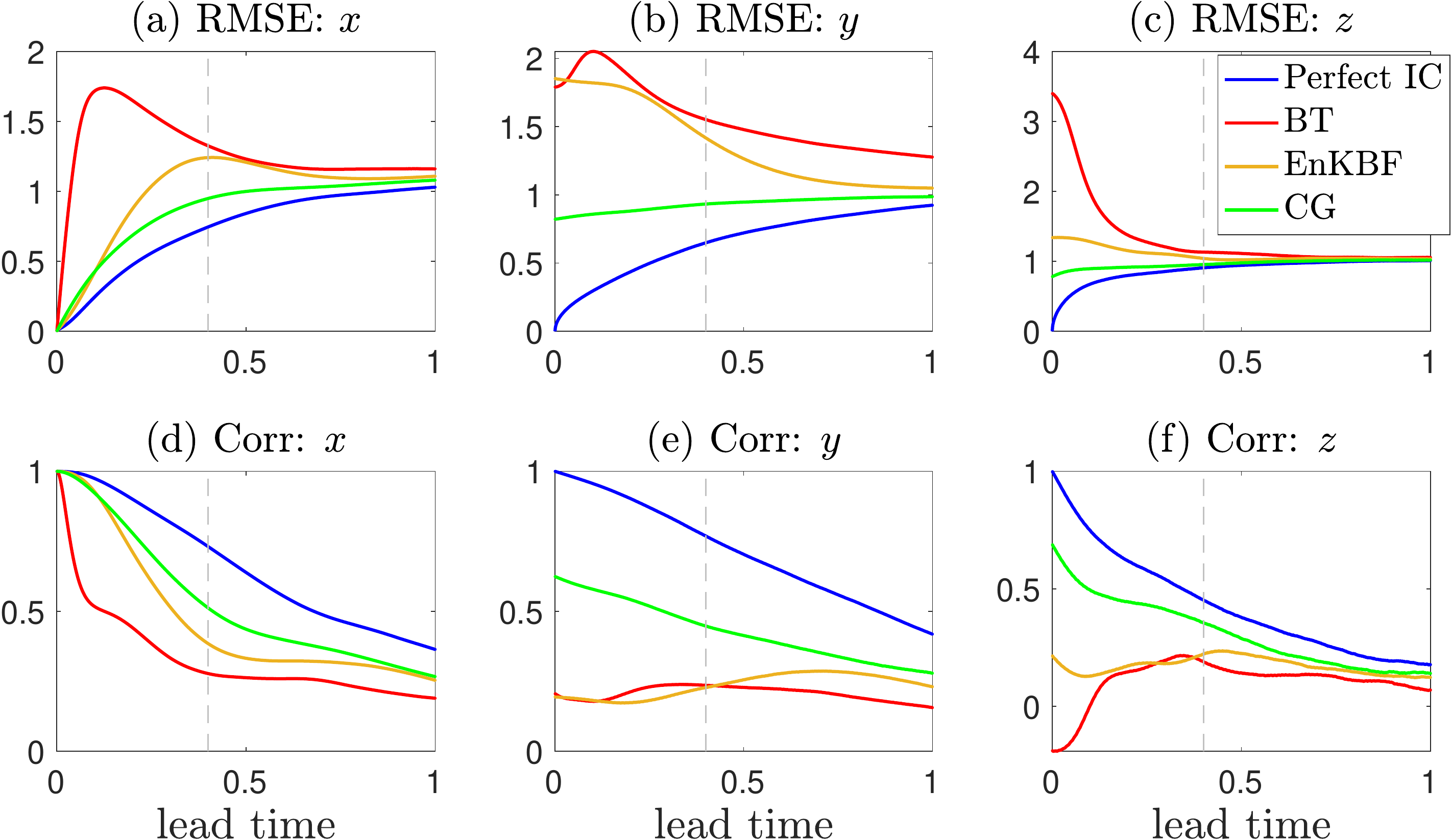}
  \caption{Panels (a)-(c): RMSE of the forecast skills for Regime II when the ICs are chosen either to be the perfect values from the true signals, or are drawn from multivariate Gaussian distributions in which the mean and covariance are estimated respectively from BT, CG and EnKBF; the RMSE are normalized by dividing by the standard deviation of the respective true signal. Panels (d)-(f): the correlation coefficients of the forecast skills. The vertical dashed line corresponds to lead time $\tau = 0.4$ for which the corresponding forecasted ensemble mean time series are shown in Figure~\ref{fig: TBH_fcst_traj_case2}.}\label{fig: TBH_rmse_case2}
\end{figure}

\begin{figure}[tbh!]
  \centering
  \includegraphics[width=1\textwidth]{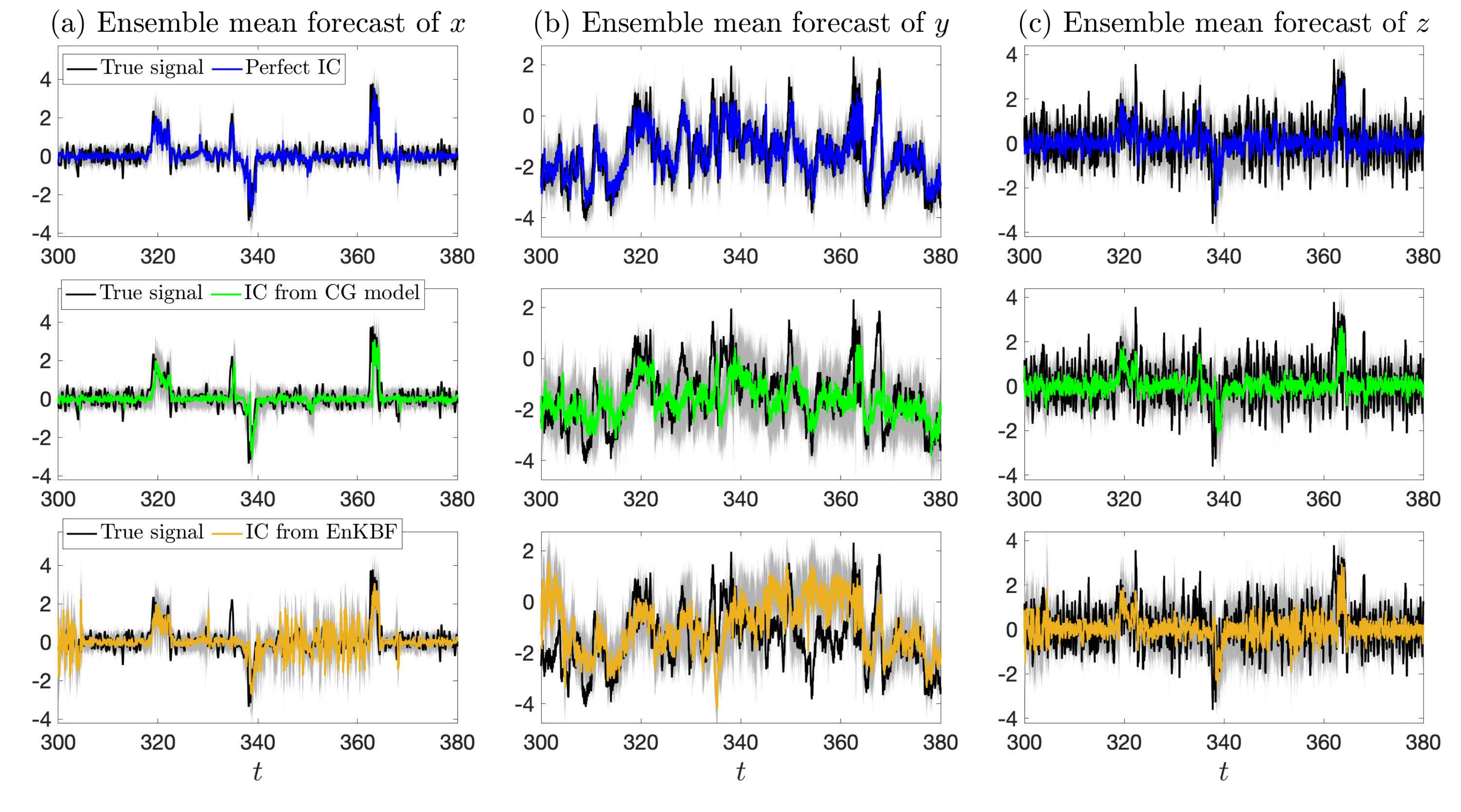}
  \caption{The forecasted ensemble mean time series for Regime II at lead time $\tau = 0.4$ time unit, when the ICs are from either the true signal (top row), or assimilated from the CG method (middle row), or the EnKBF (bottom row). The gray area on each panel marks the spread between 5 percentile and 95 percentile of the corresponding ensemble forecast. The true signals are plotted in black. For BT, large spurious oscillations appear in the forecasted ensemble mean time series for all the variables; and the corresponding results are thus not shown here.}\label{fig: TBH_fcst_traj_case2}
\end{figure}

In Figure~\ref{fig: TBH_rmse_case2}, we presented the normalized root-mean-square error (RMSE; normalized by the standard deviation of the true signal) and correlation coefficients of the forecasts for all the methods used. As is expected, the better skills of CG at the DA stage carries over to the ensemble forecast as well. BT performs the worst due to large spurious oscillations appearing in its DA stage. While the RMSE is a convenient way of ranking the performance, to provide a better visualization of the skills, we also show in Figure~\ref{fig: TBH_fcst_traj_case2} the forecasted ensemble mean trajectories at a given lead time, chosen here to be $\tau=0.4$ time unit, which corresponds roughly to one half of the decorrelation time of the $y$-variable. The results in Figure~\ref{fig: TBH_fcst_traj_case2} show that the forecast based on initial conditions provided by CG (middle row) actually performs fairly well for all the three variables compared with those when perfect initial conditions are used (top row), although the uncertainty in the $y$ variable is slightly higher for CG. The results for EnKBF are correlated with its performance at the DA stage, with significant error in the $x$ and $y$ variables over the time windows when the filtered posterior mean of $y$ deviates from the true signal as was previously shown in the middle panel of Figure~\ref{fig: TBH_DA_case2}. For BT, large spurious oscillations appear in the forecasted ensemble mean time series for all the variables, especially for $x$ and $y$. Such oscillations are inherited from those appearing in the assimilated $y$ variables, which propagate to the other two variables due to nonlinear interactions.

\section{Parameter Estimation}\label{Sec:Parameter_Estimation}

Parameter estimation is an important topic and a necessary precursor for effective state estimation, data assimilation, and prediction. Maximum likelihood estimation (MLE) and maximum a posteriori (MAP) estimation are often adopted to infer the parameter values if the observed time series for all the state variables are accessible~\cite{myung2003tutorial, kristensen2004parameter, coleman2006bayesian}. However, only partial observations are available in many complex nonlinear systems. In such a situation, data augmentation is widely used to simultaneously estimate the model parameter and recover the unobserved state variables~\cite{tanner1987calculation}. In particular, data augmentation has been extensively incorporated into the Markov Chain Monte Carlo (MCMC) algorithms for improving the Bayesian inference~\cite{wei1990monte, golightly2008bayesian, eraker2001mcmc, papaspiliopoulos2013data}. In contrast to targeting the global optimal solution based on the MCMC, many other parameter estimation algorithms seek locally optimal solutions. One of the widely used local optimal parameter estimation approaches is the expectation-maximization (EM) algorithm~\cite{ghahramani1999learning,ghahramani1996parameter,dembo1986parameter,mclachlan2007algorithm}. The EM algorithm is an iteration method that aims to find the parameter values that maximize the likelihood function and compute certain statistical expectations of the unobserved state variables in an alternating fashion. Note that, due to the local optimality property, the EM algorithm often requires fewer iterations than the MCMC algorithm. Unfortunately, both methods are computationally expensive for general complex nonlinear systems since neither the data augmentation in MCMC nor the computation of the statistical expectation in the EM can be easily obtained.

Unlike the general complex nonlinear systems, the closed analytic formulae offered by the CGNS facilitate the acceleration of the computational efficiency in parameter estimation. In particular, the conditional sampling formula~\eqref{Sampling_Main} has the potential to allow a rapid data augmentation in the MCMC method while the analytic state estimation formula~\eqref{Smoother_Main} offers an exact and accurate way to compute the statistical expectation, which is an essential component in the EM algorithm. Note that appropriately incorporating the CGNS into the MCMC may require several additional manipulations, which deserves a separate topic to study. Therefore, the focus below is on applying the CGNS to accelerate the parameter estimation utilizing the EM algorithm, where detailed mathematical justifications are more accessible. The CGNS in this entire procedure acts as a preconditioner to seek a suitable approximation of the statistical expectation of the unobserved state variables given the observational time series via the EM algorithm. With such a statistical expectation available from the CGNS, the original complex nonlinear system is only essential in computing the maximum likelihood solution at the last iteration step, leading thus to a significant speedup of the computational time required.

\subsection{Accelerating the EM algorithm with a CGNS preconditioner}\label{sec: EM_CGNS_illustration}
As before, denote by $(\mathbf{X},\mathbf{Y})$ the state variables of a complex nonlinear system of the form \eqref{Split_Form_General_Eqn}, where only a time series of $\mathbf{X}$ is observed while there is no direct observations for the state variable $\mathbf{Y}$. Let $\boldsymbol\theta$ be a vector consisting of all the model parameters. Denote by $\widehat{\mathbf{X}} = \{\mathbf{X}^0,\ldots,\mathbf{X}^j,\ldots,\mathbf{X}^J\}$ and $\widehat{\mathbf{Y}} = \{\mathbf{Y}^0,\ldots,\mathbf{Y}^j,\ldots,\mathbf{Y}^J\}$ a discrete approximation of the continuous time series of $\mathbf{X}$ and $\mathbf{Y}$, respectively, within the time interval $t\in[0,T]$, where $T=J\Delta{t}$, $\mathbf{X}^j = \mathbf{X}(t_j)$ and $\mathbf{Y}^j = \mathbf{Y}(t_j)$, with $t_j = j \Delta t$.

The goal is to seek an optimal estimation of the unknown parameters $\btheta$ by maximizing the log-likelihood function. Since only the time series of $\mathbf{X}$ is observed (denoted by the discrete sequence $\widehat{\mathbf{X}}$ above),  the log-likelihood estimate is obtained by averaging over the state variable $\mathbf{Y}$ at the corresponding instants,
\begin{equation}\label{eq: log likelihood}
\mathcal{L}(\btheta) = \log q(\widehat{\X} | \btheta) = \log \int_{\widehat{\Y}} p(\widehat{\X}, \widehat{\Y} |\btheta) \d \widehat{\Y}.
\end{equation}
Due to the unknown state variable $\mathbf{Y}$, there is in general no simple formula for a direct calculation of the log-likelihood.
To find the optimal parameter $\btheta$ that maximizes $\mathcal{L}$, the EM iteration algorithm \cite{chen2020learning,dembo1986parameter, kokkala2014expectation} operates on the following  function instead which is obtained as a lower bound of $\mathcal{L}(\btheta)$ by applying Jensen's inequality (cf.~\cite[Section 3.1]{chen2020learning})
\begin{equation}
\int_{\widehat{\Y}} Q(\widehat{\Y}) \log p(\widehat{\Y}, \widehat{\X} | \btheta) \d \widehat{\Y} - \int_{\widehat{\Y}}Q(\widehat{\Y})\log Q(\widehat{\Y})\d \widehat{\Y},
\end{equation}
where $Q(\widehat{\Y})$ is a distribution over the unobserved variable $\widehat{\Y}$. It alternates between performing an expectation (E) step to update $Q(\widehat{\Y})$ and a maximization (M) step to update $\btheta$ until the estimated parameter $\btheta$ converges.
Denote by $\btheta_k$ the updated parameters after the $k$-th iteration. The EM algorithm at step $k+1$ is the following:
\begin{itemize}
\item[\textbf{E-Step.}] Computing the conditional distribution $p(\widehat{\Y} | \widehat{\X}, \btheta_k)$ using the previously estimated parameters $\btheta_k$. In fact, the maximization in the E-Step is reached when $Q(\widehat{\Y})$ is exactly the conditional distribution of $\widehat{\Y}$ corresponding to the smoother estimates.

\item[\textbf{M-Step.}] Updating the parameters $\btheta_{k+1}$ by maximizing the cost function $\cQ$ defined by
\begin{equation}\label{eq: E-step}
  \cQ(\btheta; \btheta_k) = \int_{\widehat{\Y}} p(\widehat{\Y} | \widehat{\X}, \btheta_k) \log p(\widehat{\Y}, \widehat{\X} | \btheta) \d \widehat{\Y}.
\end{equation}
That is,
\begin{equation}\label{eq: M_step_update}
\btheta_{k+1} =\arg\max_{\btheta} \cQ(\btheta; \btheta_k).
\end{equation}
\end{itemize}

Note that in \eqref{eq: E-step}, $p(\widehat{\Y} | \widehat{\X}, \btheta_k)$ is treated as a known distribution that is solved in the E-step. Therefore, $\cQ(\btheta; \btheta_k)$ is a function of $\btheta$ only. The distribution $p(\widehat{\Y} | \widehat{\X}, \btheta_k)$ can be regarded as the weight function for computing the average (i.e., the integration) of $\log p(\widehat{\Y}, \widehat{\X} | \btheta)$.

In many situations, the M-Step usually involves solving a quadratic optimization problem, the analytic formula of which is  available. However, for general nonlinear systems, the conditional distribution $p(\widehat{\Y} | \widehat{\X}, \btheta_k)$ in the E-Step is extremely difficult to solve. Note that such a conditional distribution is precisely the smoother estimate of the complex nonlinear system. Particle methods can be applied. Yet, repeatedly using these particle methods through the iteration procedure can be computationally expensive, and careful tuning is required, especially for systems with large dimensions.

To overcome this most significant barrier of computing the conditional distribution $p(\widehat{\Y} | \widehat{\X}, \btheta_k)$ in the above EM algorithm, we exploit a suitable CGNS model as a preconditioner to accelerate this calculation in the E-Step. To that end, denote by $\btheta^M$ the collection of the parameters in the chosen CGNS model. In stark contrast to $p(\widehat{\Y} | \widehat{\X}, \btheta_k)$, we can now use the closed analytic formulae in~\eqref{Smoother_Main} to efficiently compute the conditional distribution $p^M(\widehat{\Y} | \widehat{\X}, \btheta^M_k)$ for the CGNS at each EM iteration.\footnote{We should have used $p^M(\widehat{\Y}^M | \widehat{\X}, \btheta^M_k)$ to denote the conditional distribution for the CGNS given the observation $\widehat{\X}$, where $\widehat{\Y}^M$ is the analogue of $\widehat{\Y}$ for the CGNS. We used $p^M(\widehat{\Y} | \widehat{\X}, \btheta^M_k)$ instead to avoid excessive notations. \label{Y_notation_footnote}} When $\btheta^M_k$ is converged after a sufficient number of EM iterations based on the CGNS, to obtain an approximation of the optimal $\btheta$ for the original system, we simply perform another EM cycle, but this time using $\log p(\widehat{\Y}, \widehat{\X} | \btheta)$ instead of $\log p^M(\widehat{\Y}, \widehat{\X} | \btheta^M)$ in the M-Step. The above procedure is summarized in Algorithm~\ref{Alg:EM}; and we refer to Appendix~\ref{sec: SI_EM} for further technical details. Throughout this section, a quantity with the superscript $M$ indicates that it is a quantity related to the approximate CGNS model instead of the original model (see also footnote \ref{Y_notation_footnote} with the exception of the notation $\widehat{\Y}$). Note also that $\btheta^M$ can be different from $\btheta$, since the CGNS can involve new parameters not appearing in the original system and vice versa.
\begin{algorithm2e*}[h]\label{Alg:EM}
	\caption{EM algorithm for nonlinear systems with CGNS as a preconditioner}
	\label{alg: EM_CGNS}
	Start with a given realization of the (partial) observations $\widehat\X$ for \eqref{Split_Form_General_Eqn}\;
	Propose an approximate CGNS model of the form \eqref{CGNS}\;
  Assign an initial guess for the CGNS model's parameters $\btheta^M_0$\;
  \For{$k = 1: K$}{
    E-step: compute the conditional distribution $p^{M}(\widehat{\Y} | \widehat{\X}, \btheta^M_{k - 1})$ using the previously estimated parameters $\btheta^M_{k - 1}$\;
    M-step: update the parameters $\btheta^M_{k}$ with $\btheta^M_{k} =\arg\max_{\btheta^M} \widetilde{\cQ}(\btheta^M;\btheta^M_{k - 1})$, where $\widetilde{\cQ}(\btheta^M;\btheta^M_{k - 1}) = \int_{\widehat{\Y}} p^{M}(\widehat{\Y} | \widehat{\X}, \btheta^M_{k - 1}) \log p^{M}(\widehat{\Y}, \widehat{\X} | \btheta^M) \d \widehat{\Y}$.
  }
  E-step: compute the conditional distribution $p^{M}(\widehat{\Y} | \widehat{\X}, \btheta^M_{K})$\;
  M-step: compute the optimal parameters $\btheta_{\text{opt}}$ for \eqref{Split_Form_General_Eqn} with $\btheta_{\text{opt}} =\arg\max_{\btheta} \cQ(\btheta;\btheta^M_{K})$, where $\cQ(\btheta;\btheta^M_{K}) = \int_{\widehat{\Y}} p^{M}(\widehat{\Y} | \widehat{\X}, \btheta^M_{K}) \log p(\widehat{\Y}, \widehat{\X} | \btheta) \d \widehat{\Y}$.
\end{algorithm2e*}

Note that $\btheta_{\text{opt}}$ obtained from Algorithm~\ref{Alg:EM} should be viewed as an approximation of the true optimal parameters for the original nonlinear system. Its quality relies obviously on the quality of the CGNS in approximating the true dynamics. It is also worth mentioning that additional (physics) constraints can be naturally incorporated into the CGNS for parameter estimation, while still preserving the closed analytic formulae in the corresponding EM algorithm; and under certain conditions, a block decomposition of the conditional covariance can be devised to further reduce the computational efforts when high dimensional systems are considered~\cite{chen2020learning}. We illustrate now the approach on the classical two-layer Lorenz 1996 model.

\subsection{A multiscale turbulent test model} \label{sec:L96}

In this subsection, a two-layer inhomogeneous Lorenz model is utilized to demonstrate that a suitable CGNS can serve as both a preconditioner and a surrogate model in parameter estimation to speed up the computation. First, we show that a simple approximate model that belongs to CGNS can accelerate the EM algorithm as a preconditioner by following Algorithm~\ref{Alg:EM}. Then, we show that this approximate model itself can be used as a surrogate model for prediction once the involved parameters, $\btheta^M$, are optimized through the EM cycles in lines 4-6 of Algorithm~\ref{Alg:EM}.

\subsubsection{The perfect model}
The two-layer Lorenz 96 (L96) model \cite{lorenz1996predictability} is a conceptual representation of geophysical turbulence that is commonly used as a testbed for various stochastic parameterization and dimension reduction techniques \cite{chorin2015discrete,fatkullin2004computational,arnold2013stochastic,crommelin2008subgrid,vissio2018proof,majda2009mathematical,grooms2014stochastic,grooms2014stochastic2,majda2014new}. The model mimics a coarse discretization of atmospheric flow on a latitude circle. It supports complex wave-like and chaotic behavior, and the two-layer structure schematically depicts the interactions between small-scale fluctuations and large-scale motions. The stochastic version of the model subject to additive noise forcing reads

\begin{subequations}\label{eq: L96_model_inhom}
	\begin{align}
		\frac{d u_{i}}{d t} &= \left(- u_{i - 1}\left(u_{i - 2}-u_{i + 1}\right)-u_{i} + f - \frac{h c_i}{J} \sum^{J}_{j=1} v_{i, j} \right) + \sigma_{u_i} \dot W_{u_i}, \quad i = 1, \dots, I, \label{eq: L96_model_inhom_a} \\
		\frac{d v_{i, j}}{d t} &= \left(-b c_i v_{i, j + 1}\left(v_{i, j + 2} - v_{i, j - 1}\right) - c_i v_{i, j} + \frac{h c_i}{J} u_{i} \right) + \sigma_{v_{i, j}} \dot W_{v_{i, j}}, \quad j=1, \dots, J, \label{eq: L96_model_inhom_b}
\end{align}
\end{subequations}
where $I$ denotes the total number of large-scale variables, $J$ the number of small-scale variables corresponding to each large-scale variable, $f$, $h$, $c_i$, $b$, $\sigma_{u_i}$ and $\sigma_{v_{i,j}}$ are given scalar parameters while $\dot W_{u_i}$ and $\dot W_{v_{i, j}}$ are white noise. The large-scale variables~$u_i$ are periodic in $i$ with $u_{i + I} = u_{i - I} = u_{i}$. The corresponding small-scale variables~$v_{i, j}$ are periodic in $i$ with $v_{i + I, j} = v_{i - I, j} = v_{i, j}$ and satisfy the following cyclic conditions in $j$: $v_{i, j + J} = v_{i + 1, j}$, and $v_{i, j - J} = v_{i - 1, j}$.

The model discussed here uses variables $u_i$ to describe large-scale or slow movements which are resolved; small scales or rapid fluctuations represented by $v_{i,j}$ are often unresolved ones. The coupling of fast and slow variables is regulated by the parameter $h$. The parameter $f$ controls the magnitude of external large-scale forcing, while $b$ determines the amplitude of nonlinear interactions between the fast variables. The parameter $c_i$ specifies how quickly the fast variables are damped in comparison to the slow variables. We take $I=40$, corresponding to a discretization of the latitude circle into a total of 40 sites/sectors, and choose $J = 4$ small-scale variables associated with each $u_i$. The constant forcing is set to be $f = 4$, which makes the system chaotic for the parameter regime chosen here. The parameters $h$, $c_i$, and $b$ are chosen in such a way that the small-scale variables have a comparatively significant impact on the large-scale ones. In other words, the perfect model only has a weak scale separation. The reason that we consider such a weak scale separation is that it better mimics the real atmosphere with chaotic/turbulent behavior. The parameter $c_i$ varies across the spatial sites, which aims to mimic the fact that the coupling across the variables above the ocean is weaker than that above the land since the latter usually have stronger friction or dissipation. In this sense, the model is inhomogeneous. Finally, additional stochastic noise is added to the system, representing the contribution of the variables that are not explicitly modeled. The noise also interacts with the deterministic part via nonlinear terms, introducing additional complexity that mimics nature. To summarize, the parameters used in the perfect model~\eqref{eq: L96_model_inhom} are as follows,
\begin{equation}\label{eq: L96_model_para}
	\begin{aligned}
    I = 40, \quad J = 4, \quad h = 2, \quad c_i = 2 + 0.7\cos(2 \pi i / I), \quad b = 2, \quad  f = 4, \quad \\
    \sigma_{u_i} = \sigma_u = 0.2, \quad \sigma_{v_{i, j}} =  \sigma_v = 1. \qquad \qquad \qquad \qquad
  \end{aligned}
\end{equation}

\subsubsection{The approximate model}
Since in general the perfect model is not always fully known, or it is too complicated to be used in practice, it is essential to develop a simple and computationally tractable approximate model, which is nevertheless able to capture the key nonlinear feedback from the unobserved variables ($v_{i, j}$ here) to the observed variables ($u_i$ here). As was discussed in Section~\ref{Sec:Model_Development}, stochastic parameterization is widely used in describing chaotic signals~\cite{chen2016filtering}, which replaces the nonlinear eddy terms by quasilinear stochastic processes on formally infinite embedded domains where the stochastic processes are Gaussian conditional to the large scale mean flow. In addition, physics-constraint is adopted in designing  approximate models, which also includes the effects from the large-scale $u_i$ to the small-scale variables~$v_{i, j}$. Therefore, such approximate models can potentially be used as surrogate models of the perfect model.

The approximate model that we introduce for \eqref{eq: L96_model_inhom} is as follows
\begin{subequations}\label{eq: L96_SP}
  \begin{align}
  \frac{d u_{i}}{d t} &= \left(- u_{i - 1}\left(u_{i - 2}-u_{i + 1}\right)-u_{i} + \hat{f}_i - \hat{a}_i \sum^{J}_{j=1} v_{i, j} \right) + \hat{\sigma}_{u_i} \dot W_{u_i}, \quad i = 1, \dots, I, \label{eq: L96_SP_1} \\
  \frac{d v_{i, j}}{d t} &= - \hat{d}_{i, j}v_{i, j} + \hat{v}_{i, j} +  \hat{c}_i u_i + \hat{\sigma}_{v_{i, j}} \dot W_{v_{i, j}}, \quad j=1, \dots, J,~\label{eq: L96_SP_2}
  \end{align}
\end{subequations}
where $\hat{f}_i$,  $\hat{a}_i$, $\hat{\sigma}_{u_i}$, $\hat{d}_{i, j}$, $\hat{v}_{i, j}$, $\hat{c}_i$, and $\hat{\sigma}_{v_{i, j}}$ are unknown constants.
Compared with the original system \eqref{eq: L96_model_inhom}, the main simplification here is in the small-scale equations, where we have replaced the small-scale nonlinear interaction and linear damping terms $-b c_i v_{i, j + 1}(v_{i, j + 2} - v_{i, j - 1}) - c_i v_{i,j}$ in \eqref{eq: L96_model_inhom_b} by a simple linear term $- \hat{d}_{i,j} v_{i,j} + \hat{v}_{i, j}$. The equations for $u_i$ are essentially the same as before, although $\hat{f}_i$ is allowed now to vary from site to site. In the numerical experiments below, we will enforce $\hat{a}_i = \hat{c}_i$ for the consistency of the dynamical property with the original system \eqref{eq: L96_model_inhom}.

One desirable feature of this approximate model is that the direct coupling of the state variables only involves $u_i$ and the corresponding $v_{i, j}$ for each fixed $i$. This is different from the original system \eqref{eq: L96_model_inhom} where $u_i$ can have direct interactions with $v_{i+1,1}$ and $v_{i-1,J}$ due to the cyclic boundary conditions of the small scale $v_{i, j}$ over $j$. Such a property allows to use a block decomposition of the covariance matrix of the smoother estimate during both E-step and M-step~\cite{chen2020learning}. The entire state space for all the variables $\{u_i, v_{i,j}\; \vert \;  i=1, \dots, I, j = 1, \dots J\}$ can be decomposed into $I$ subspaces, where each subspace can be dealt with in parallel.

\subsubsection{Setup of the numerical simulations}
The true signal is obtained by integrating the inhomogeneous L96 model~\eqref{eq: L96_model_inhom} using the Euler-Maruyama scheme with the parameters given by \eqref{eq: L96_model_para}, a uniform time step size $\delta t = 2 \times 10^{-3}$, and zero initial condition. The same time step size and initial condition are adopted for the approximate model \eqref{eq: L96_SP} as well as the identified model, where the latter takes the same form as \eqref{eq: L96_model_inhom} but with estimated parameters. The true signals of $u_i$ with 100 time units are used as the observations for parameter estimation while longer data with 1000 time units are used to compute their statistics. This latter length is also adopted when computing the statistics of the variables in the approximate and the identified models. The total number of the EM loops with CGNS is fixed to be 200. 

\subsubsection{CGNS as a preconditioner for identifying parameters in the perfect model}
\begin{figure}[tbh!]
  \centering
  \includegraphics[width=1\textwidth]{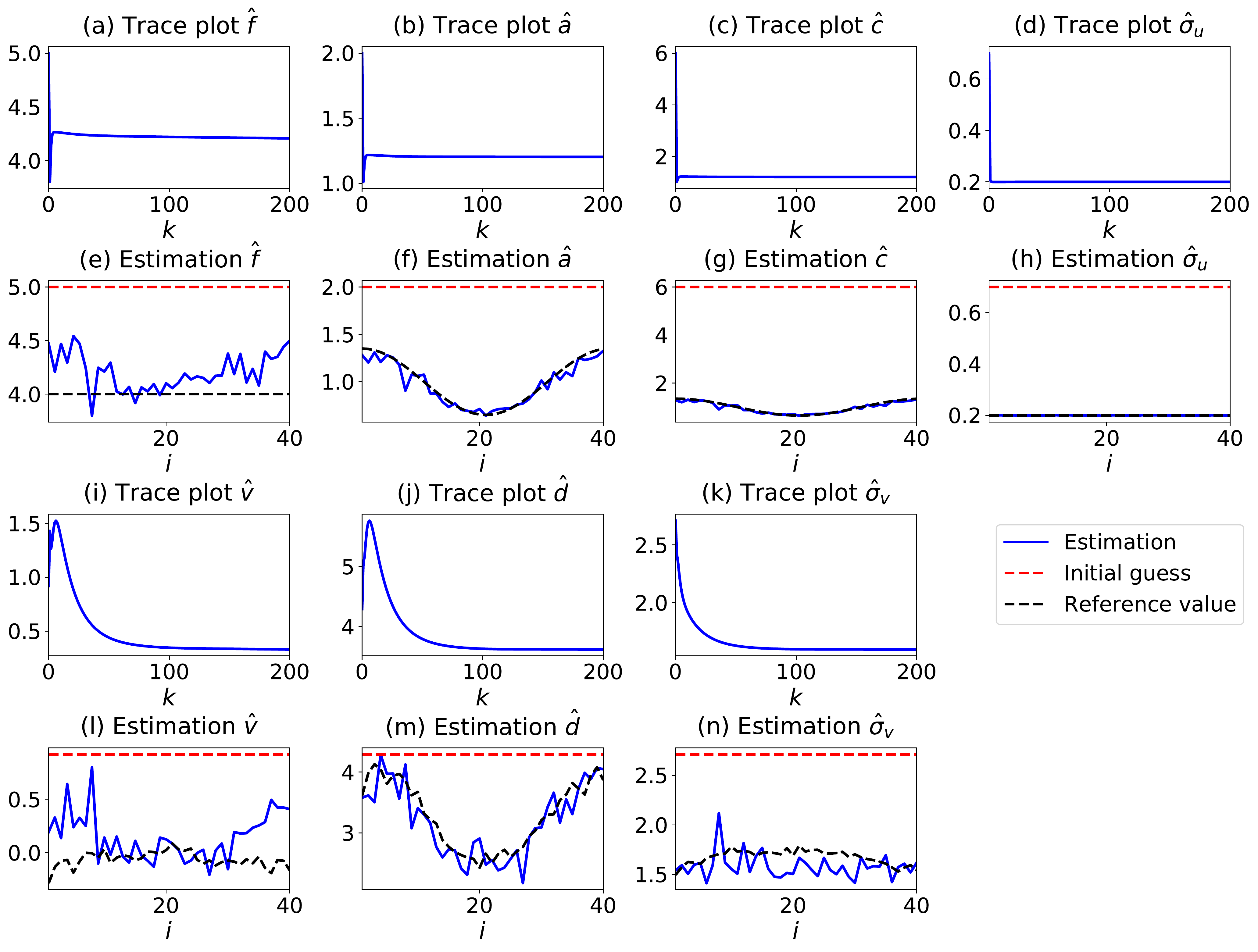}
  \caption{Learning the parameters of the approximate L96 model~\eqref{eq: L96_SP}. The panels (a)--(d) and (i)--(k) in the 1st and 3rd rows show the trace plots of the estimated parameters of the EM algorithm~\ref{Alg:EM} at site $i=2$. The panels (e)--(h) and (l)--(n) in the 2nd and 4th rows show the final estimated parameters~$\btheta^M_{K}$ (blue) after the $K$-th step of the EM algorithm with $K$ taken here to be 200; also shown are the initial guesses of the parameters (red) as well as the reference parameter values (back) chosen as follows: $\hat{f}_i$, $\hat{a}_i$ ($=hc_i /J$), $\hat{\sigma}_{u_i}$, and $\hat{c}_i$ are set to be the corresponding true values used in the perfect model, while $\hat{d}_{i,j} = \hat{d}_{i}$, $\hat{v}_{i, j} = \hat{v}_{i}$, and $\hat{\sigma}_{v_{i, j}} = \hat{\sigma}_{v_{i}}$ are calibrated by the true statistics according to \eqref{eq: calculate_params_in_L96_SP}.}\label{fig: trace_plot_L96}
\end{figure}

We first discuss the results of applying Algorithm~\ref{Alg:EM} to estimate the parameters, $\btheta^M$, in the approximate model~\eqref{eq: L96_SP}. Figure~\ref{fig: trace_plot_L96} shows that the EM algorithm with the approximate model \eqref{eq: L96_SP} provides a very accurate approximation of the parameters corresponding to the Gaussian fit of $v_{i,j}$ in the true signal. The quantities shown in Figure~\ref{fig: trace_plot_L96} are the trace plots at the site $i = 2$ corresponding to the variables $u_2$ and $v_{2,j}$ (in the first and the third rows) and the final estimation of the parameters in~\eqref{eq: L96_SP} (in the second and the fourth rows). The trace plots at other spatial locations have similar behavior; the parameters involved in the $u_i$-equations all converge quickly (within 10 iteration steps), and those involved in the $v_{i,j}$-equations converge at a relatively slower speed, but all stabilized after about 100 iterations. The black curves are shown as a reference, where $\hat{f}_i$, $\hat{a}_i$ ($=hc_i /J$), $\hat{\sigma}_{u_i}$, and $\hat{c}_i$ are taken to be those in the perfect model and $\hat{d}_{i,j} = \hat{d}_{i}$, $\hat{v}_{i, j} = \hat{v}_{i}$, and $\hat{\sigma}_{v_{i, j}} = \hat{\sigma}_{v_{i}}$ are calibrated by the true statistics, i.e., the mean, the variance, and the decorrelation time, in the perfect model of $v_{i, j}$ averaged over $j$. More precisely, denoting by $\mu_i$, $r_i$, and $\tau_i$ these averaged mean, variance, and decorrelation time in the perfect model of $v_{i,j}$. the parameters $\hat{v}_i$, $\hat{d}_i$, and $\hat{\sigma}_{v_{i}}$ are then obtained via
\begin{equation}\label{eq: calculate_params_in_L96_SP}
  \tau_i = \frac{1}{\hat{d}_i}, \qquad r_i = \frac{\hat{\sigma}_{v_{i}}^2}{2\hat{d}_i},  \qquad \mu_i = \frac{\hat{v}_i + h c_i/J \bar{u}_i}{\hat{d}_i} ,
\end{equation}
where $\bar{u}_i$ is the mean value of $u_i$.

\begin{figure}[tbh!]
  \centering
  \includegraphics[width=1\textwidth]{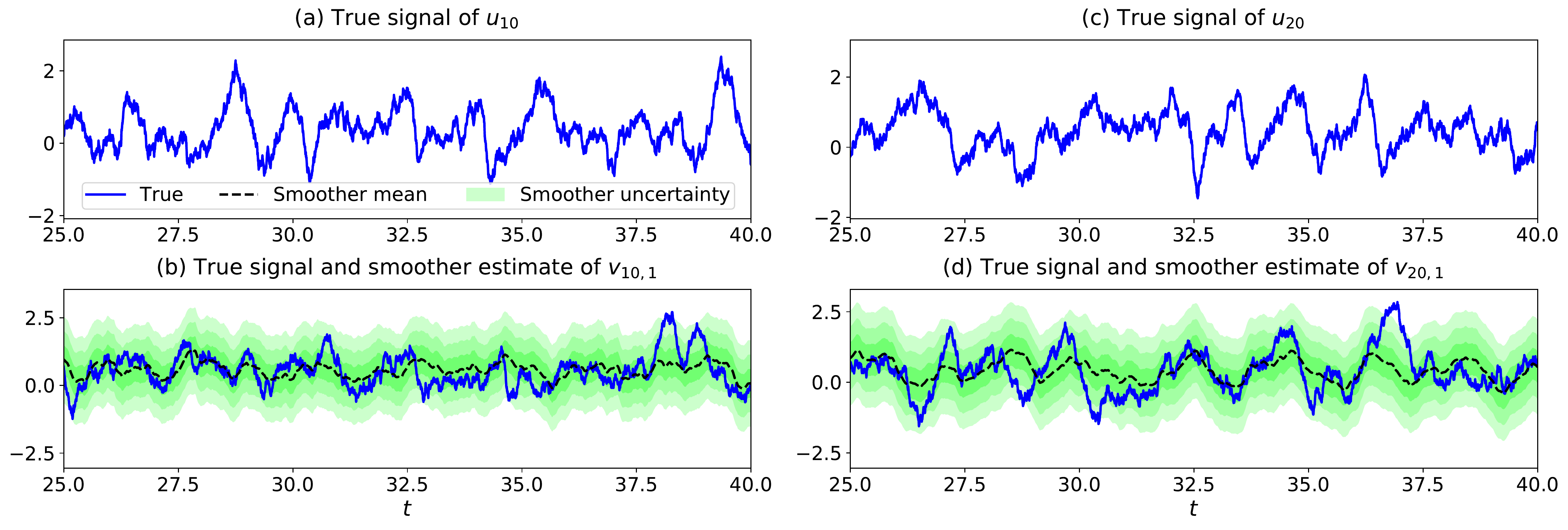}
  \caption{Smoother estimate of the approximate model in the $K$-th iteration in Algorithm~\ref{Alg:EM}, with $K=200$. The blue curves: true signals; the black dashed curves: the smoother mean time series of the hidden variable; the light, moderate, and dark shading areas show the one, two, and three standard derivations (STDs) of the uncertainty in the smoother estimate. Panel (a): true signal of $u_{10}$; Panel (b): true signal of $v_{10, 1}$ with one, two, and three, standard derivations of the uncertainty. Panel (c)--(d): Same as Panel (a)--(b) but for $u_{20}$ and $v_{20, 1}$.}\label{fig: smoother_dist_L96}
\end{figure}

The conditional distribution of the approximate model~\eqref{eq: L96_SP} with the estimated parameters is shown in Figure~\ref{fig: smoother_dist_L96}. Panel (a) and (c) show the true signal of two large-scale modes $u_{10}$ and $u_{20}$, and Panel (c) and (d) show the true signal, and smoother estimate of two hidden modes. The smoother mean is given by the dashed-black curves, and one, two, and three standard derivations of the uncertainty are shown by the light, moderate, and dark shading areas. The shading areas cover most of the true signal, which indicates appropriate amount of uncertainties are obtained by combining the true observations of the large scale variables $u_i$ and the optimized approximate model. In fact, characterizing appropriate amount of uncertainty plays an important role in the EM algorithm when the hidden process contains large uncertainty.  If the uncertainty is totally ignored in computing the optimization in the M-Step, for example replacing $p(\widehat{\Y} | \widehat{\X}, \btheta^M_k)$ by its mean state, then the estimated parameters can be very biased. In fact, the solution of the estimation even blows up in the test model used here.

\begin{figure}[tbh!]
  \centering
  \includegraphics[width=0.8\textwidth]{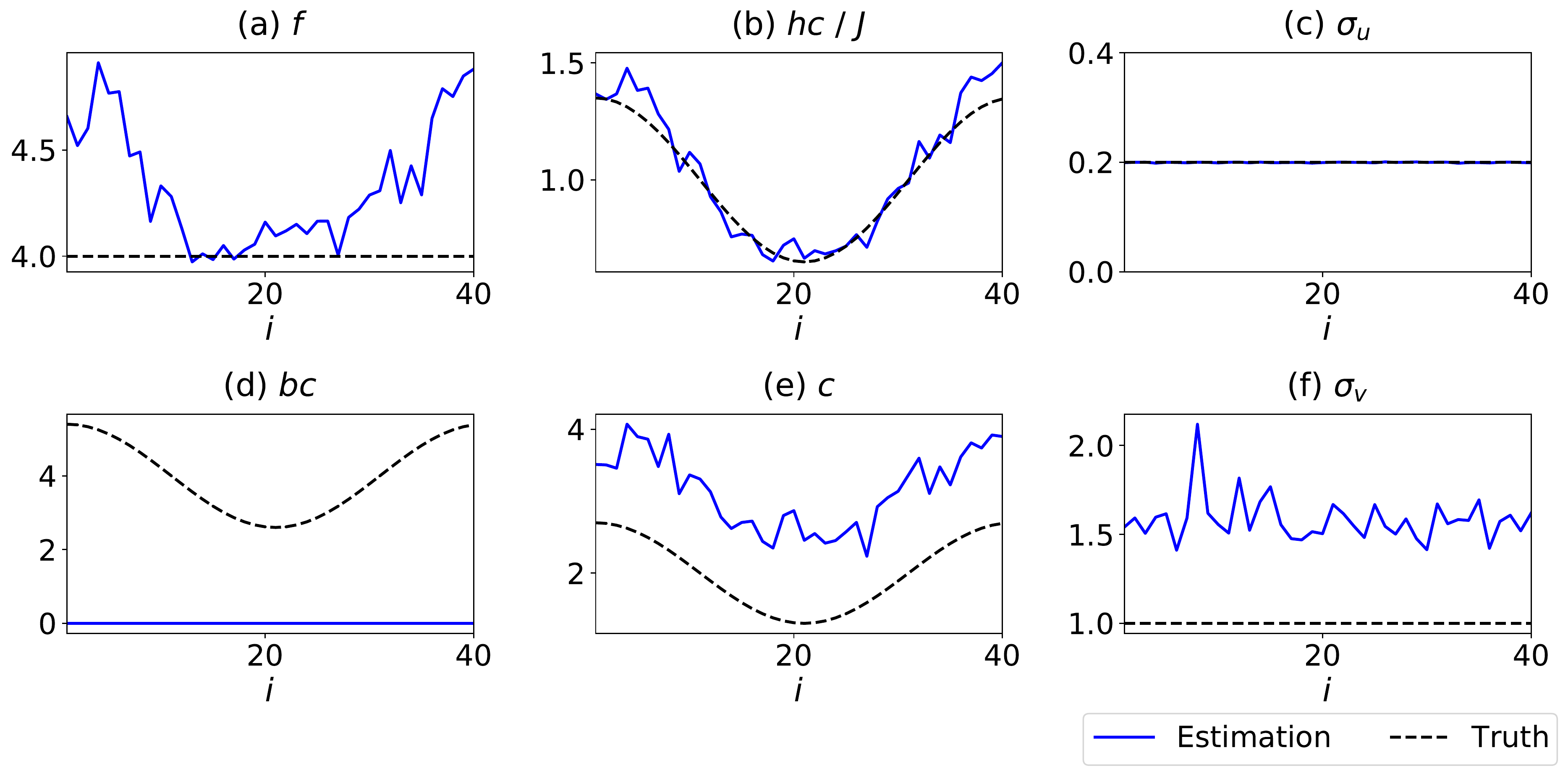}
  \caption{The estimated parameter $\btheta_{\text{opt}}$ using Algorithm~\ref{alg: EM_CGNS} for the perfect L96 model~\eqref{eq: L96_model_inhom}. Here, $\btheta_{\text{opt}}$ is a vector consisting of the parameters $f$, $hc_i/J$, $\sigma_u$, $bc_i$, $c_i$, and $\sigma_v$ in \eqref{eq: L96_model_inhom}. The true values given by \eqref{eq: L96_model_para} are marked by the dashed black curves, and the estimated parameters are shown by the blue curves.}\label{fig: estimated_params_L96_final}
\end{figure}

Figure~\ref{fig: estimated_params_L96_final} shows the estimated parameters $\btheta_{\text{opt}}$ for the perfect model, where the CGNS \eqref{eq: L96_SP} is utilized as a preconditioner following Algorithm \ref{Alg:EM}. Due to the intrinsic model error of the approximate CGNS model where the hidden variables are fully decomposed, in the sense that the correlation between the small scales corresponding to different large scales are omitted, the estimated value for the $b c_i$ term is zero. However, other parameters, for example $f$ and $c_i$, are adjusted accordingly. It is shown in Figure~\ref{fig: Hovmoller_L96} and Figure~\ref{fig: recovered_traj_L96} that the identified full model with the estimated parameters $\btheta_{\text{opt}}$ produces dynamics that resembles the truth to a remarkable extent. Indeed, Figure~\ref{fig: Hovmoller_L96} shows that the hovmoller plot of the large scale variables from the identified model (Panel (b)) is almost the same as that from the true signal (Panel (a)) over the given time window. Time series comparison as well as PDF and ACF comparisons are also shown in the first four row of Figure~\ref{fig: recovered_traj_L96}. Both the dynamical properties and the statistics are recovered with high accuracy. Of course, the time series from the identified model should not be expected to follow the true time series at all time instants since the original model is placed in a regime with chaotic dynamics.

\subsubsection{CGNS as a surrogate model}
\begin{figure}[tbh!]
  \centering
  \includegraphics[width=0.8\textwidth]{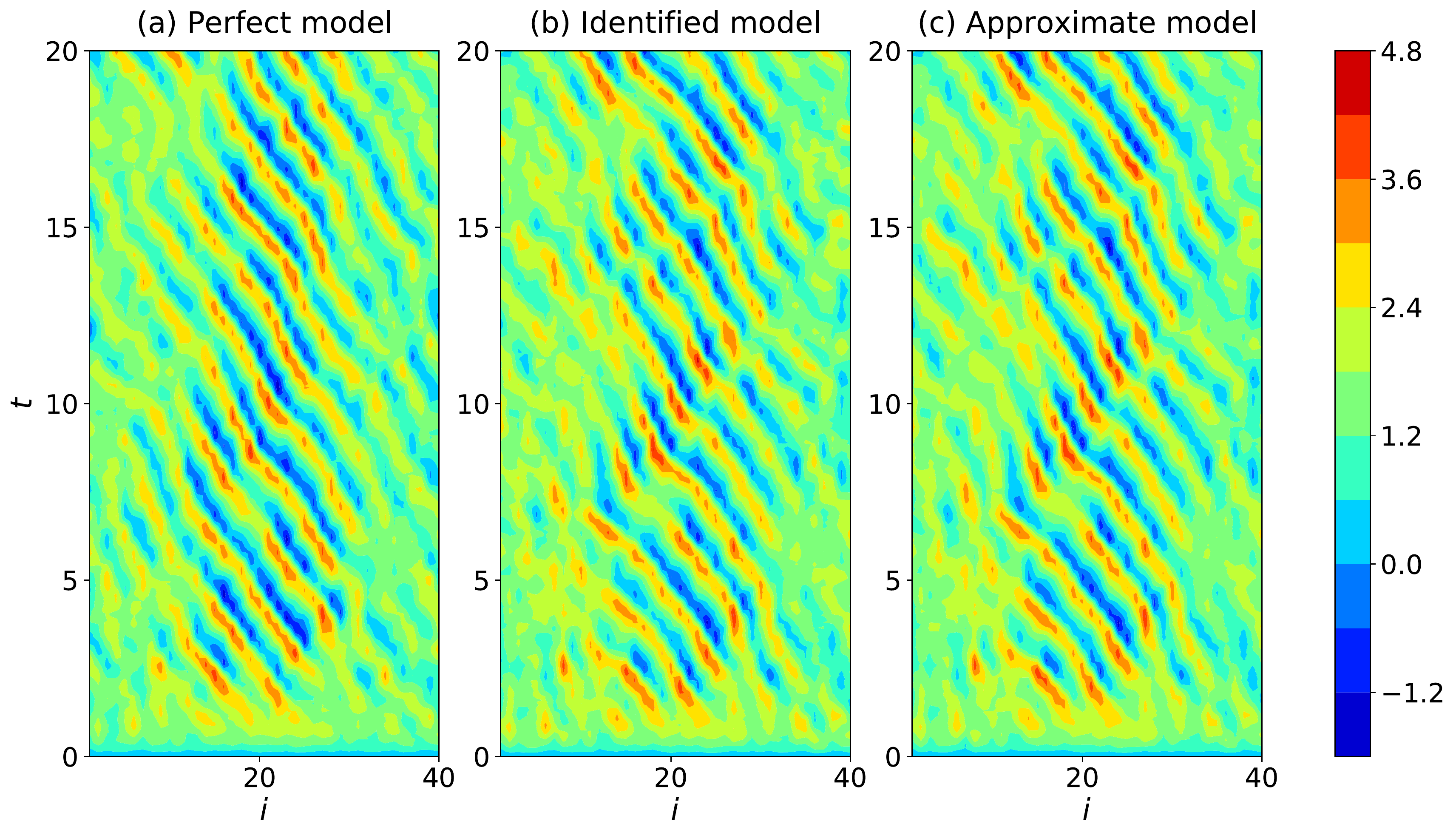}
  \caption{Hovmoller diagram of the large scale variables $u_i$ from different models. Panel (a): perfect model~\eqref{eq: L96_model_inhom} with parameters~\eqref{eq: L96_model_para}; Panel (b): perfect model~\eqref{eq: L96_model_inhom} with estimated parameters~$\btheta_{\text{opt}}$ after line 8 using Algorithm~\ref{alg: EM_CGNS}; Panel (c): approximate model~\eqref{eq: L96_SP} with parameters~$\btheta^M_K$ with $K = 200$.}\label{fig: Hovmoller_L96}
\end{figure}

\begin{figure}[tbh!]
  \centering
  \includegraphics[width=1\textwidth]{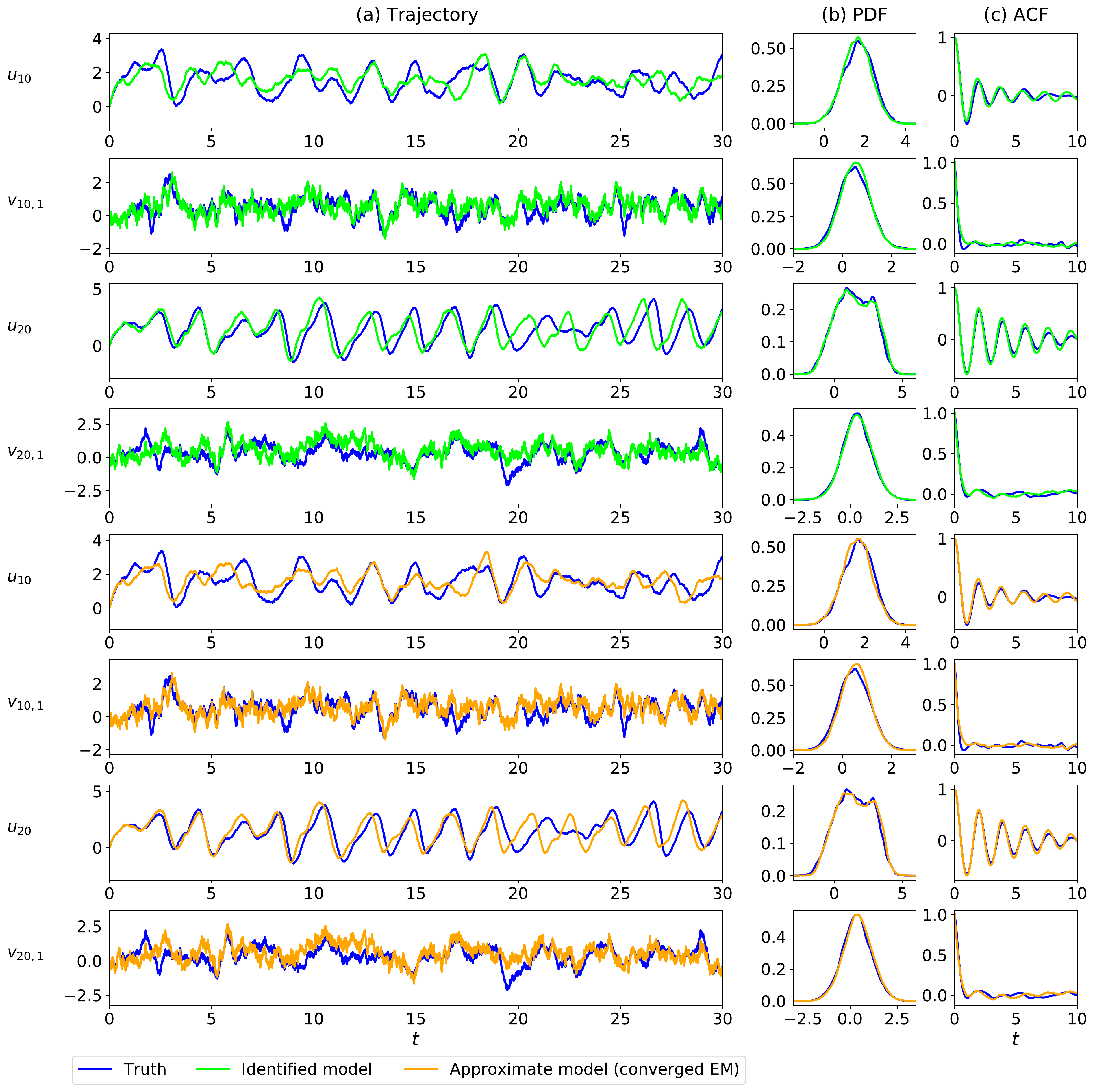}
  \caption{Comparison of time series and statistics obtained from the full two-layer L96 model~\eqref{eq: L96_model_inhom} with the true parameters~\eqref{eq: L96_model_para} (blue), the identified model~\eqref{eq: L96_model_inhom} with the estimated parameters~$\btheta_{\text{opt}}$ from Algorithm~\ref{Alg:EM} (green), and the approximate CGNS model~\eqref{eq: L96_SP} with the estimated parameters~$\btheta^M_{K}$ (orange) after the $K$-th iteration in Algorithm~\ref{alg: EM_CGNS} where $K=200$. Panel (a)--(c): trajectory, PDF, and ACF, respectively.}\label{fig: recovered_traj_L96}
\end{figure}

Finally, we mention that the approximate model~\eqref{eq: L96_SP} with the estimated parameters itself can be exploited as a surrogate model of the perfect system, which can be applied for ensemble forecast and other tasks. Indeed, \eqref{eq: L96_SP} with the optimally estimated parameters $\btheta^M_K$ from Algorithm~\ref{Alg:EM} recovers the dynamical properties of the true model to an extent that is almost the same as the full model with the identified parameters $\btheta_{\text{opt}}$ as shown in Figure~\ref{fig: Hovmoller_L96} and the last four rows in Figure~\ref{fig: recovered_traj_L96}.

\section{Predicting the Statistical Response}\label{Sec:Response}

Yet another important topic in studying complex nonlinear systems is to predict the model response to the perturbation of the external forcing. Developing efficient and accurate approaches to study such a key issue facilitates understanding  model sensitivity, regime-switching behavior, and nonlinear interactions across different scales. Resolving this problem using advanced mathematical tools also has significant practical implications, such as coping with the climate change scenario. Due to various uncertainties from the internal instability and external forcing, a probabilistic description is more suitable for characterizing the complex turbulent systems~\eqref{eq:abs_formu}.

However, there exist several challenges in predicting the statistical response of complex nonlinear systems. First, solving the high-dimensional Fokker-Planck equation is the prerequisite of obtaining the model statistics, which, however, often suffers from the curse of dimensionality. Second, since simulating the perfect system is not always computationally feasible in practice, the predicted model response may become inaccurate when an approximate model is utilized. It is then important to take advantage of a suitable combination of partial observations with the approximate model to mitigate error in calculating the model statistics. Third, for general nonlinear systems, computing the statistical response in terms of the perturbations with different strengths and categories requires repeatedly solving the Fokker-Planck equation. As a consequence, even with a relatively fast solver of the Fokker-Planck equation, the total computational cost can still remain significant.

The advantage of the exact and statistically accurate solver of the equilibrium PDF~\eqref{Joint_X_Y} makes the CGNS  a natural framework for the development of approximate models for efficiently computing the model statistics and the associated response. One important feature in finding the PDF of the CGNS based on the formula~\eqref{Joint_X_Y} is that the available partially observed time series $\mathbf{X}$ is used to compute the conditional distribution of $\mathbf{Y}$ given $\mathbf{X}$. As a consequence, the model error in the PDF associated with the direct simulation of the approximate model is mitigated with the help of observations in computing the PDF based on~\eqref{Joint_X_Y}. In other words, the resulting PDF from~\eqref{Joint_X_Y} is in general closer to the truth than that computed purely based on the approximate model without taking into account any input information from observations.

What remains is the third challenge mentioned above. To overcome such a difficulty, the linear statistical response is utilized as an approximate method for computing the exact statistical response. The linear response only requires a linearization of the statistical equation while there is no linearization involved in the original nonlinear dynamics. Therefore, the nonlinear features of the underlying dynamics is preserved. In addition, the linear response can be computed utilizing the fluctuation-dissipation theorem (FDT) \cite{majda2005information}, which involves only a single PDF for computing different linear responses. Note that such a PDF is the equilibrium distribution of the unperturbed state,  which can be efficiently solved using the formula~\eqref{Joint_X_Y_Equilibrium}.

\subsection{Computing the linear statistical response via the fluctuation-dissipation theorem (FDT)}
Consider the general complex nonlinear systems in~\eqref{eq:abs_formu}. Defining $\G(\mathbf{u},t) = \left(L+D\right)\u+B\left(\u,\u\right)+\F(t)$, the model can be written in a concise form
\begin{equation}\label{eq: FDT_general}
  \frac{\d\u}{\d t} = \mathbf{G}(\u,t) + \boldsymbol{\sigma}(\u,t)\dot{\mathbf{W}}.
\end{equation}
The equilibrium statistics of some functional $A(\u)$ associated with~\eqref{eq: FDT_general} is formulated as
\begin{equation}\label{eq: FDT_Stat}
  \langle A(\u)\rangle = \int A(\u)\peq(\u)\d\u,
\end{equation}
where $\peq(\u)$ is the equilibrium PDF of $\u$ in~\eqref{eq: FDT_general}.

Now consider the dynamics in~\eqref{eq: FDT_general} by a small time separable external forcing perturbation $\delta \mathbf{w}(\u) f(t)$, where $\delta$ is a small scalar $\mathbf{w}$ is a general nonlinear function of $\mathbf{u}$. The perturbed system reads
\begin{equation}\label{eq: FDT_general_perturb}
  \frac{\d\u}{\d t} = \mathbf{G}(\u,t) + \delta \mathbf{w}(\u) f(t) + \boldsymbol{\sigma}(\u,t)\dot{\mathbf{W}}.
\end{equation}
The FDT states that if $\delta$ is small enough, then the leading-order correction to the statistics in~\eqref{eq: FDT_Stat} is given by
\begin{equation}\label{eq: FDT_Stat2}
  \delta\langle A(\u)\rangle(t) = \delta \int_0^t\mathbf{R}(t-s) f(s)\d s,
\end{equation}
where $\mathbf{R}(t)$ is the linear response operator, which is calculated through correlation functions in the unperturbed dynamics:
\begin{equation}\label{eq: FDT_Formula}
  \mathbf{R}(t) = \langle A(\u(t))\cB(\u(0))\rangle,\qquad \cB (\u) = -\frac{\mbox{div}_\u(\mathbf{w (\u)}\peq(\u))}{\peq(\u)}.
\end{equation}
See \cite[Chapter 2]{majda2005information} for a rigorous derivation of~\eqref{eq: FDT_Stat2}--\eqref{eq: FDT_Formula}. In particular, the above procedure of computing the linear response via the FDT does not require the linearization of the underlying complex nonlinear systems. Therefore, the features of the nonlinear dynamics are preserved. Note also that if the functional $A(\u)$ in~\eqref{eq: FDT_Stat2} is given by $A(\u)=\u$, then the response computed is for the statistical mean. Likewise, $A(\u)=(\u-\bar{\u})^2$ is used for computing the response in the variance.

\subsection{Calculating the linear statistical response via the CGNS preconditioner}\label{sec: FDT_CGNS_procedure}
According to~\eqref{eq: FDT_Formula}, the calculation of the linear statistical response via the FDT requires the information of
\begin{enumerate}
  \item the equilibrium PDF $\peq(\u)$,
  \item the time series $A(\u)$, and
  \item the correct formulation of~$\cB(\u)$.
\end{enumerate}
Even if the perfect model is known, directly solving the high-dimensional Fokker-Planck equation is often not computationally affordable. Therefore, a suitable CGNS, serving as a preconditioner, is utilized to find a suitable approximation of the non-Gaussian equilibrium PDF $\peq$ in an efficient way. Specifically, the explicit formula in~\eqref{Joint_X_Y_Equilibrium} is utilized to achieve this goal. Besides, the observations  also need to be incorporated in computing the equilibrium $\peq$ (and later in recovering the hidden components in $A(\u)$) to reduce model error from the approximate model free run. We denote one realization of the observational variable by $\X^{\obs}$, the posterior distribution (filter and smoother) given $\X^{\obs}$ by~$p^{M |\obs}$ where the superscript $M$ indicates that an approximate model is used in computing the filter distribution~\eqref{CGNS_Stat}, the explicit formula~\eqref{Joint_X_Y_Equilibrium} becomes
\begin{equation}\label{Joint_X_Y_Equilibrium_Imperfect}
  \peq^{M | \obs}(\mathbf{X},\mathbf{Y}) =\lim_{J\to\infty} \frac{1}{J}\sum_{j=1}^J \Big(K_\mathbf{H}(\mathbf{X}-\mathbf{X}^{\obs}(t_j)) p^{M | \obs}(\mathbf{Y}|\mathbf{X}^{\obs}(s\leq t_j))\Big),
\end{equation}
where $\peq^{M | \obs}$ is an efficient and effective approximation of the true equilibrium PDF~$\peq$.

Next, in the presence of partial observations, the conditional sampling formula~\eqref{Sampling_Main} is exploited to calculate the unobserved component of the time series in $A(\u)$. Note that the approximate model is used to compute the filter distribution~\eqref{CGNS_Stat}, the smoother distribution formulae~\eqref{Smoother_Main} and the conditional sampling formula~\eqref{Sampling_Main}. Note that the partial observations are involved in computing both $\peq(\u)$ and $A(\u)$, which aims to  mitigate the model error in the approximate model in both the equilibrium PDF and time series of the unobserved components. In addition, with $\peq(\u)$ and $A(\mathbf{u})$ obtained from the CGNS preconditioner, the original nonlinear system structure is utilized to form $\cB(\u)$ to compute the linear response $\R(t)$. The entire procedure of the FDT via CGNS is given in Algorithm~\ref{alg: FDT_CGNS}.

\begin{algorithm2e}[H]
	\caption{FDT with the CGNS preconditioner}
	\label{alg: FDT_CGNS}
	Start with a given realization of the observations $\X^{\obs}$\;
	Propose an approximate model that belongs to CGNS~\eqref{CGNS}\;
  Compute the filter posterior distribution $p^{M | \obs}(\Y | \X^{\obs}(s \le t_i))$ via~\eqref{CGNS_Stat} \;
	Form the equilibrium $\peq^{M|\obs}$ via equation~\eqref{Joint_X_Y_Equilibrium_Imperfect}\;
  Compute the smoother posterior distribution~$p^{M | \obs}(\Y(t)|\X^{\obs}(s), s\in[0,T])$ from~\eqref{Smoother_Main}\;
  Sample one realization of the hidden time series via~\eqref{Sampling_Main} that is used to approximate the unobserved component of $A(\mathbf{u})$\;
  Compute the response operator $\mathbf{R}(t)$ via~\eqref{eq: FDT_Formula} and compute the linear response~via~\eqref{eq: FDT_Stat2}.
\end{algorithm2e}

\subsection{A 4D stochastic climate model}

This section utilizes a four-mode stochastic model with key features of atmospheric low-frequency variability to show how to use CGNS as a preconditioner incorporated with partial observations to calculate the linear statistical response.
\subsubsection{The perfect model}
The stochastic climate model is designed in such a way that it involves many of the major dynamical properties of comprehensive global circulation models (GCMs) but with only four degree of freedom~\cite{majda2008applied,majda2005information,majda1999models,majda2001mathematical}. The model reads as follows
\begin{subequations}\label{eq: 4d_model_CG}
  \begin{align}
    \frac{\d x_1}{\d t} &= \left(-x_2(L_{12} + a_1 x_1 + a_2x_2) - d_1 x_1 + F_1 + L_{13}y_1 + b_{123}x_2 y_1 \right) + \sigma_1 \dot{W}_{x_1},\\
    \frac{\d x_2}{\d t} &= \left(+ x_1(L_{12} + a_1 x_1 + a_2x_2) - d_2 x_2 + F_2 + L_{24}y_2 + b_{213}x_1 y_1 \right) + \sigma_2 \dot{W}_{x_2}, \\
    \frac{\d y_1}{\d t} &= \left(- L_{13}x_1 + b_{312} x_1 x_2 + F_3 - \gamma_1 y_1 \right) + \sigma_3 \dot{W}_{y_1}, \\
    \frac{\d y_2}{\d t} &= \left(- L_{24}x_2 + F_4 - \gamma_2 y_2\right) + \sigma_4  \dot{W}_{y_2},
  \end{align}
\end{subequations}
where $b_{123} + b_{213} + b_{312} = 0$. Consistent with many geophysical flow models, the model has energy-conserving quadratic nonlinear terms, a linear operator, and external forcing terms. The linear operator contains two parts: one is a skew-symmetric component formally related to the Coriolis effect and topographic Rossby wave propagation; the other is a negative definite symmetric portion conceptually analogous to dissipative processes such as surface drag and radiative damping. The coupling in different variables is through both linear and nonlinear terms, where the nonlinear coupling through $b_{ijk}$ produces multiplicative noise effects. In fact, the strategies described in Section~\ref{Subsec:averaging} are applied to $y_1$ and $y_2$ that introduce the stochastic noise and damping terms. The variables $x_1$ and $x_2$ can be regarded as the climate variables and $y_1$ and $y_2$ represents the weather variables.
The parameters used to generate the true dynamics are as follows
\begin{equation}\label{eq: 4d_model_para}
\begin{aligned}
& d_1 = 1, \quad d_2 = 0.4, \quad \gamma_1 = 0.5, \quad \gamma_2 = 0.5, \quad L_{12} = 1, \quad L_{13} = 0.5,  \quad L_{24} = 0.5, \\
  & a_1 = 2, \quad a_2 = 1, \quad b_{123} = 1.5, \quad b_{213}=1.5, \\
& \sigma_1 = 0.5, \quad \sigma_2 = 2, \quad \sigma_3 = 0.5, \quad \sigma_4 = 1, \quad F_1 = F_2 = F_3 = F_4 = 0. \\
\end{aligned}
\end{equation}
One realization of the true signal is shown in black color in Figure~\ref{fig: 4d_climate_traj}. Both climate variable $x_1$ and weather variable $y_1$ have intermittent behavior with non-Gaussian PDFs. Note that this stochastic model is CGNS with $\X=(x_1, x_2)^{\top}$ and $\Y=(y_1, y_2)^{\top}$. We use a CGNS as the perfect model such that the FDT based on the perfect model can be computed in an accurate fashion, which can be served as a reference solution.

\subsubsection{The approximate model}
In practice, running the entire perfect model is prohibitively costly. As a result, simpler or reduced models are commonly utilized in computing the responses. Linear stochastic models are widely used as approximate models for the unresolved variables~\cite{delsole2005predictability}. Therefore, the hidden processes are replaced by two linear Gaussian equations, the parameters of which are calibrated by the true equilibrium statistics, i.e., the mean, the variance, and the decorrelation time. Besides, the parameters in the observed processes are assumed to be the same as those in the perfect model~\eqref{eq: 4d_model_CG}. The approximate model reads,
\begin{subequations}\label{eq: 4d_model_CG_SP}
  \begin{align}
    \frac{\d x_1}{\d t} &= \left(-x_2(L_{12} + a_1 x_1 + a_2x_2) - d_1 x_1 + F_1 + L_{13}y_1 + b_{123}x_2 y_1 \right) + \sigma_1 \dot{W}_{x_1},\\
    \frac{\d x_2}{\d t} &= \left(+ x_1(L_{12} + a_1 x_1 + a_2x_2) - d_2 x_2 + F_2 + L_{24}y_2 + b_{213}x_1 y_1 \right) + \sigma_2 \dot{W}_{x_2}, \\
    \frac{\d y_1}{\d t} &= - \hat{\gamma}_3 \left(y_1 - \hat{y}_1 \right) + \hat{\sigma}_3  \dot{W}_{y_1}, \\
    \frac{\d y_2}{\d t} &=- \hat{\gamma}_4 \left(y_2 - \hat{y}_2 \right) + \hat{\sigma}_4  \dot{W}_{y_1},
  \end{align}
\end{subequations}
which belongs to the CGNS. Note that despite the simplicity of utilizing linear Gaussian models to approximate the hidden processes, one major issue in \eqref{eq: 4d_model_CG_SP} is that the physics constraint is no longer satisfied in~\eqref{eq: 4d_model_CG_SP}. Therefore, the model in \eqref{eq: 4d_model_CG_SP} can contain large errors for a long-term simulation.

\subsubsection{Calculating linear response via FDT}
In this example, we aim at calculating the linear response to the perturbation of the external forcing and linear interaction parameters. Specifically, the following two perturbation cases are considered:
\begin{compactitem}
  \item[] \textbf{Case 1}: Perturbing parameters of forcing in the observed processes, i.e., $F_1^{\delta} = F_2^{\delta} = 0.3$, $\delta \w(\u) f(t) = (0.3, 0.3, 0, 0)^{\top}$;
  \item[] \textbf{Case 2}: Perturbing parameters in linear interaction terms, i.e., $L_{13}^{\delta} = L_{24}^{\delta} = 0.1$, $\delta \w(\u) f(t) = (0.1 y_1, 0.1 y_2, -0.1 x_1, -0.1 x_2)^{\top}$.
\end{compactitem}
Note that in the second case, the parameters appear in both the observed and hidden processes.
We compare the linear response in the following models:
\begin{compactitem}
  \item[]  \textbf{Perfect FDT (or perfect model)}: The equilibrium PDF~$\peq(\u)$, the time series $A(\u)$ and the formulation of $\cB(\u)$ are all from the perfect model~\eqref{eq: 4d_model_CG}.
  \item[]  \textbf{Imperfect FDT (or imperfect model)}: The equilibrium PDF~$\peq(\u)$, the time series $A(\u)$ and the formulation of $\cB(\u)$ are all from the approximate model~\eqref{eq: 4d_model_CG_SP}.
  \item[]  \textbf{Concatenate FDT  (or concatenate model)}:  The equilibrium PDF~$\peq(\u)$, the time series $A(\u)$ are from the simple concatenation of the observations, i.e., the trajectories of observed variables $x_1$ and $x_2$ from the perfect model~\eqref{eq: 4d_model_CG}, and the trajectories of the hidden variables $y_1$ and $y_2$ from the approximate model~\eqref{eq: 4d_model_CG_SP} free-run. The formulation of $\cB(\u)$ is from the approximate model~\eqref{eq: 4d_model_CG_SP}.
  \item[]  \textbf{FDT with preconditioner (or preconditioner model)}: The equilibrium PDF~$\peq(\u)$ is calculated via~\eqref{Joint_X_Y_Equilibrium} where true observations and the CGNS approximate model \eqref{eq: 4d_model_CG_SP} are utilized in computing~$K_\mathbf{H}(\X-\X^{\obs}(t_i))$ and~$p^{M | \obs}(\Y|\X(s\leq t_i))$. The time series $A(\u)$ are generated by equation~\eqref{Sampling_Main} with the CGNS approximate model and the true observations. The formulation of $\cB(\u)$ is from the perfect model.
\end{compactitem}
The details of the $\cB(\u)$'s forms from the perfect model and approximate model can be found in Appendix~\ref{sec: Bu_form}.
In this experiment, the true signal is obtained using Euler-Maruyama scheme with a uniform time step~$\delta t = 5 \times 10 ^{-3}$ and the total length is 1000 time units.

\begin{figure}[tbh!]
  \centering
  \includegraphics[width=1\textwidth]{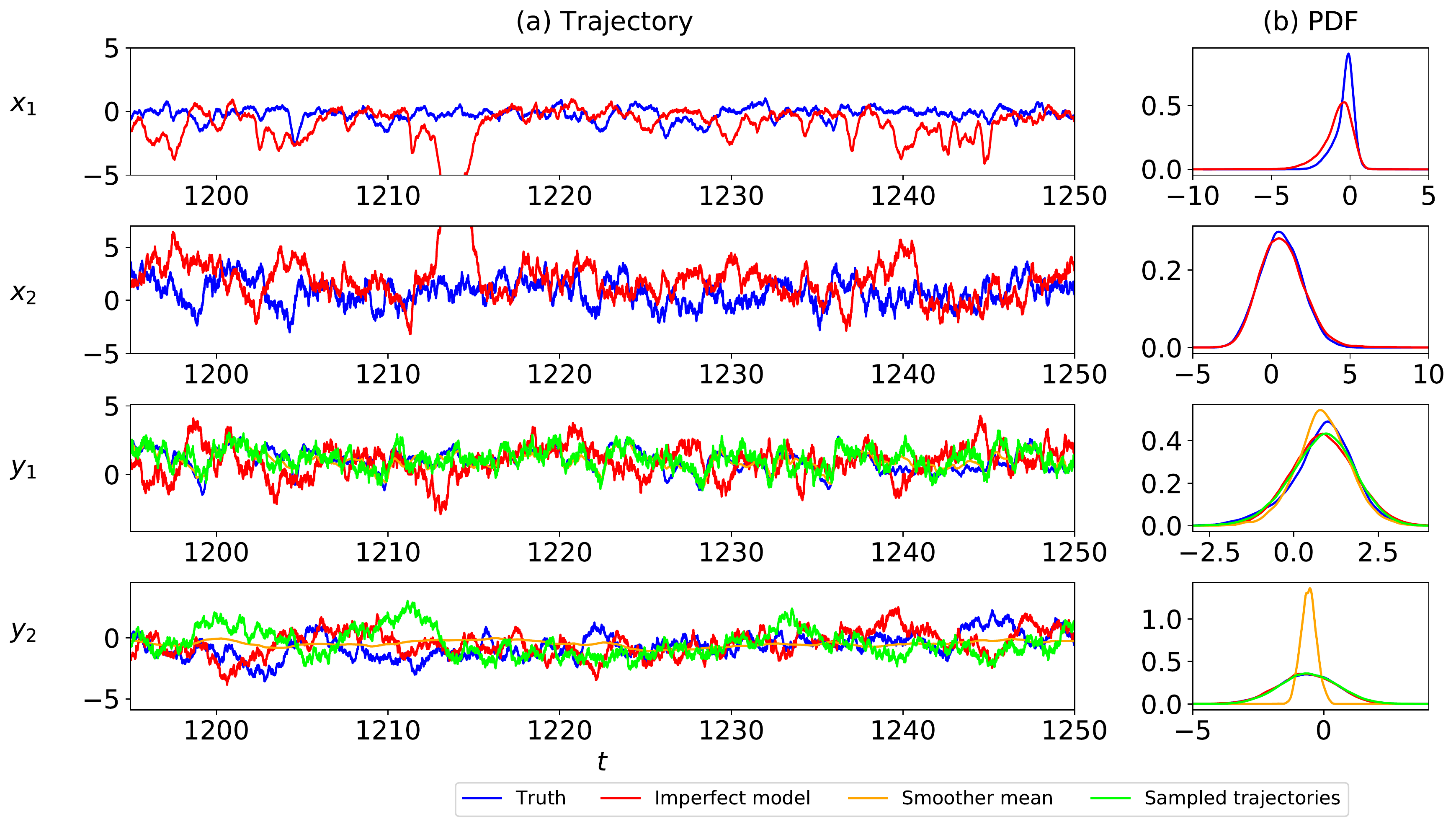}
  \caption{Comparison of the trajectories from the perfect 4D stochastic climate model~\eqref{eq: 4d_model_CG}, approximate model~\eqref{eq: 4d_model_CG_SP}, the smoother mean time series, and the sampled trajectories based on the approximate model. Blue curves: perfect model trajectories; red curves: approximate model trajectories; orange curves: the smoother mean time series; green curves: sampled trajectories. Panel (a): trajectories; Panel (b): PDFs.}\label{fig: 4d_climate_traj}
\end{figure}

\begin{figure}[tbh!]
  \centering
  \includegraphics[width=1\textwidth]{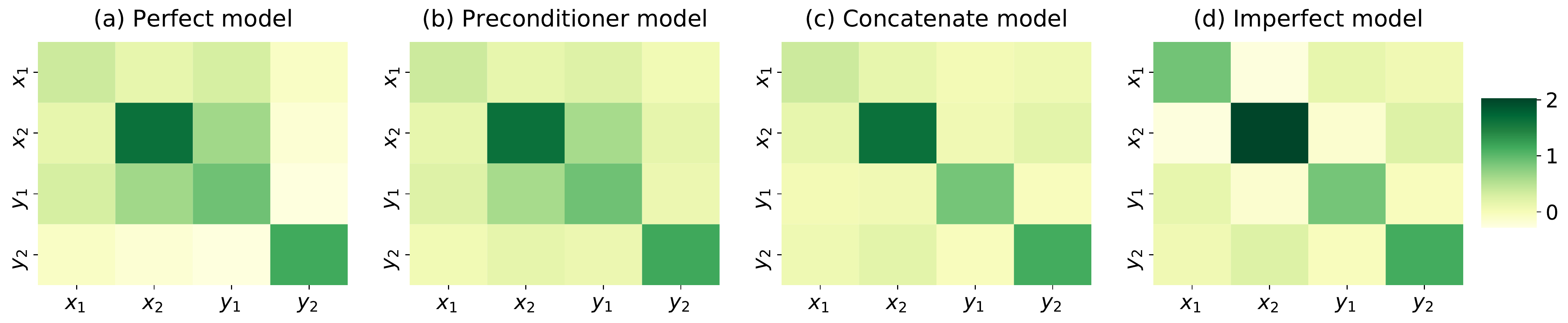}
  \caption{Comparison of covariance matrices in different scenarios. Panel (a): the perfect model where the equilibrium PDF is from the perfect model~\eqref{eq: 4d_model_CG}; Panel (b): preconditioner model where equilibrium PDF is from~\eqref{CGNS_Stat} where true observations and approximate model are utilized in computing~$K_\mathbf{H}(\X-\X^{\obs}(t_i))$ and~$p^{M | \obs}(\Y|\X(s\leq t_i))$; Panel (c): concatenate model where the equilibrium PDF~$\peq$ is from the simple concatenation of the observations, i.e., the trajectories of observed variables $x_1$ and $x_2$ from the perfect model~\eqref{eq: 4d_model_CG} and the trajectories of the hidden variables $y_1$ and $y_2$ from the approximate model~\eqref{eq: 4d_model_CG_SP} free-run; Panel (d): the imperfect model where the equilibrium PDF is from the approximate model~\eqref{eq: 4d_model_CG_SP}.}\label{fig: 4d_climate_cov}
\end{figure}

Before discussing the linear response using CGNS as a preconditioner, we start by showing some pre-requisites, the equilibrium covariance, and time series $A(\u)$ from different models. The trajectories from the approximate model~\eqref{eq: 4d_model_CG_SP} free-run (red color) are shown in Figure~\ref{fig: 4d_climate_traj}. Note that due to the violation of the energy-conserving constraint of the nonlinear terms in~\eqref{eq: 4d_model_CG_SP}, the amplitude of $x_1$ and $x_2$ from the approximate model is much larger than the one from the perfect model~\eqref{eq: 4d_model_CG}. In addition, the PDF of $y_1$ is Gaussian by design, which is also different from the skewed PDF as in the perfect model. The model error can also be found in the comparison of the equilibrium covariance matrices of the perfect model (Panel (a)) and the approximate model (Panel (d)) in Figure~\ref{fig: 4d_climate_cov}. Due to large error caused by the simple parameterization of the hidden processes, one may consider concatenating the observations (from the perfect model) and the trajectories of the hidden variables $y_1$ and $y_2$ from the approximate model~\eqref{eq: 4d_model_CG_SP} free-run to approximate the required equilibrium PDF~$\peq$ and the unobserved time series. However, it is expected that the correlation between the observations from the perfect model and trajectories from the approximate model free-run is neglected. See Panel (c) of Figure~\ref{fig: 4d_climate_cov}.

In contrast, in light of the desired structures of CGNS, the partial observations can be incorporated with both approximating the equilibrium PDF and recovering the unobserved times series. Therefore, the model error is significantly mitigated, and the correlation between the observed and unobserved variables is reserved.  The green curves in Figure~\ref{fig: 4d_climate_traj} show one sampled trajectories of $y_1$ and $y_2$ using conditional sampling formula~\eqref{Sampling_Main} from the approximate model~\eqref{eq: 4d_model_CG_SP} and observations. The overall dynamics of the recovered sampled trajectories are very similar to those in the perfect model. In addition, the skewed PDF of~$y_1$ can be found in the green curve but the appropriate model free run only brings Gaussian statistics by design. More importantly, the correlation between the sampled trajectories and the observations is consistent with that in the perfect model as shown in Panel (b) of Figure~\ref{fig: 4d_climate_cov}. A final remark is that the smoother (or filter) mean time series is widely used as a surrogate of the true signals, which, however, underestimate the uncertainty as shown in the orange color in Figure~\ref{fig: 4d_climate_traj}. Therefore, conditional sampling is essential in approximating the time series of $A(\mathbf{u})$.

\begin{figure}[tbh!]
  \centering
  \includegraphics[width=1\textwidth]{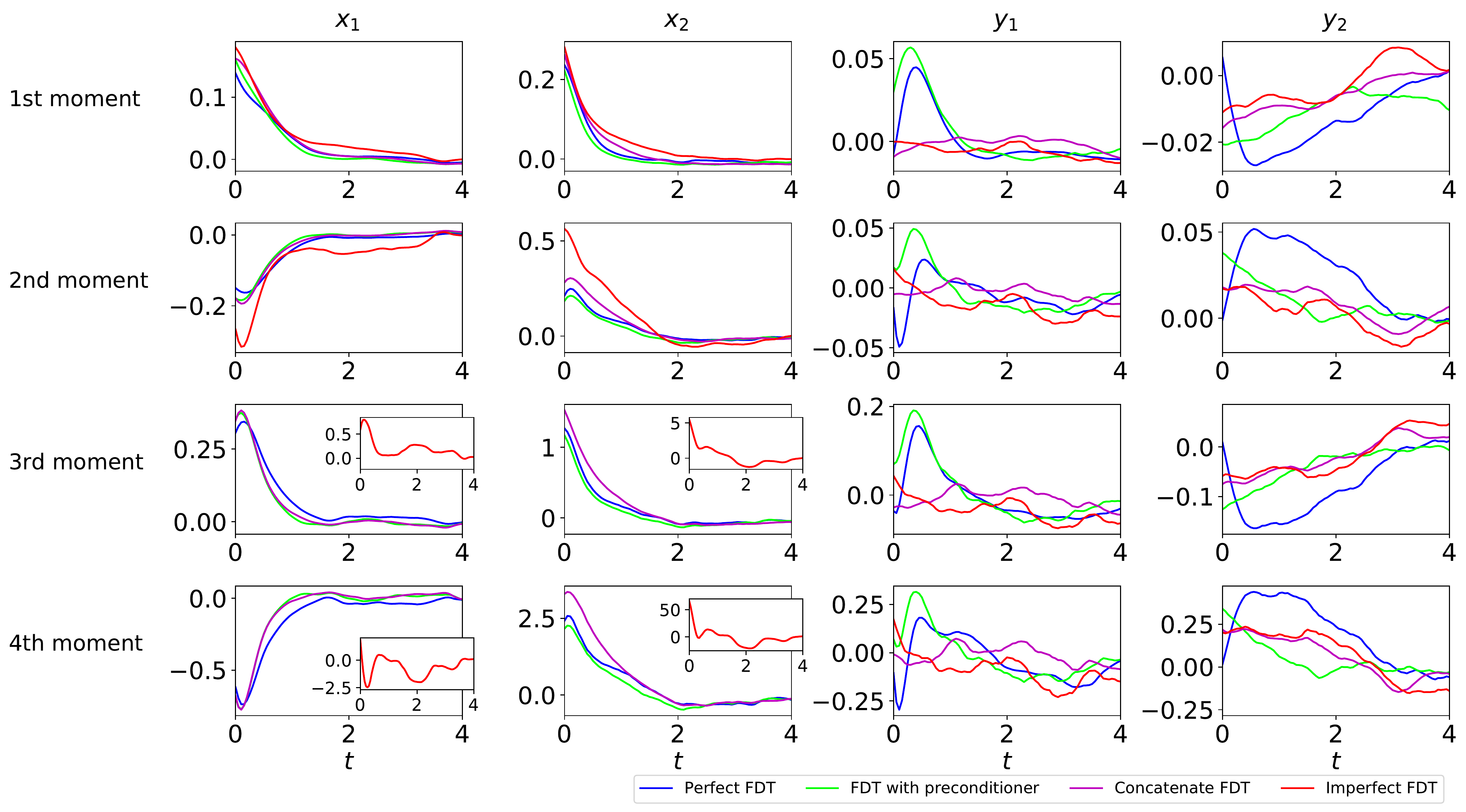}
  \caption{Response operator $\R(t)$ in~\eqref{eq: FDT_Formula} for the response of the first four moments of the 4D climate model~\eqref{eq: 4d_model_CG} when perturbing parameters in the observed processes with $F_1^{\delta} = F_2^{\delta} = 0.3$. In each panel, the blue and red curves show the linear response from the perfect model~\eqref{eq: 4d_model_CG} and the free-run of the approximate model~\eqref{eq: 4d_model_CG_SP}, respectively. The magenta curves show linear response from simply concatenating the trajectories of observed variables $x_1$ and $x_2$ from the perfect model~\eqref{eq: 4d_model_CG} and the trajectories of the hidden variables $y_1$ and $y_2$ from the approximate model~\eqref{eq: 4d_model_CG_SP}. The green curves show the linear response from the procedure discussed in Section~\ref{sec: FDT_CGNS_procedure}. Panel (a)--(d): $x_1$, $x_2$, $y_1$, and $y_2$, respectively. }\label{fig: 4d_climate_response_obs}
\end{figure}

\begin{figure}[tbh!]
  \centering
  \includegraphics[width=1\textwidth]{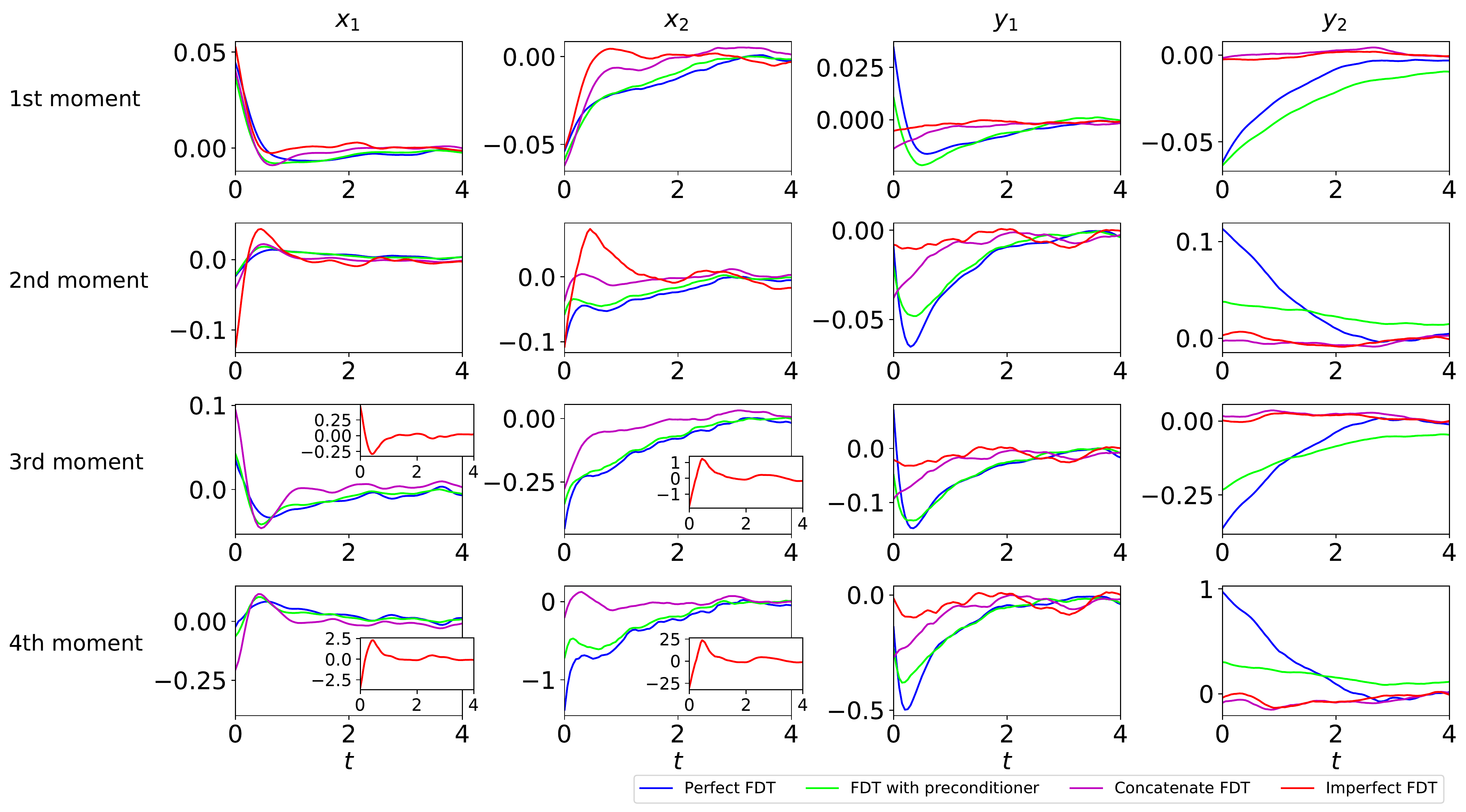}
  \caption{Similar to Figure~\ref{fig: 4d_climate_response_obs} but perturbing parameters in the linear interaction terms $L_{13}^{\delta} = L_{24}^{\delta} = 0.1$. Panel (a)--(d): $x_1$, $x_2$, $y_1$, and $y_2$, respectively.}\label{fig: 4d_climate_response_hidden}
\end{figure}

Figure~\ref{fig: 4d_climate_response_obs} shows the linear response operator $R(t)$ in~\eqref{eq: FDT_Formula} for the response of the first four moments when perturbing the external forcing parameters in the observed processes. Here, the blue and red curves show the linear response from the perfect model~\eqref{eq: 4d_model_CG} and the free-run of the approximate model~\eqref{eq: 4d_model_CG_SP}, respectively. The magenta curves show the linear response from simply concatenating the trajectories of observed variables $x_1$ and $x_2$ from the perfect model~\eqref{eq: 4d_model_CG} and the trajectories of the hidden variables $y_1$ and $y_2$ from the approximate model~\eqref{eq: 4d_model_CG_SP}. The green curves show the linear response from the procedure discussed in Section~\ref{sec: FDT_CGNS_procedure}. The imperfect FDT contains huge errors in capturing higher moments for the observed variables $x_1$ and $x_2$. Concatenate FDT works well for computing the four moments of $x_1$; however, it is gradually away from the perfect FDT from lower moments to higher moments. For example, the gap between the concatenate FDT and the perfect FDT is obvious in the last row of $x_2$. This is because the strong correlation between the $x_2$ and $y_1$ (shown in Panel (a) in Figure~\ref{fig: 4d_climate_cov}) is omitted in this simple concatenation (shown in Panel (c) in Figure~\ref{fig: 4d_climate_cov}). In contrast, the response operator $R(t)$ of the two observed variables from preconditioner FDT is very close to the perfect FDT. Due to the indirect perturbation of the unresolved variables, the response is expected to be small. The preconditioner FDT can still capture the trend of the linear response operator as that using perfect FDT. Figure~\ref{fig: 4d_climate_response_hidden} shows the second perturbation case, i.e., perturbing the linear interaction parameters appearing in both the observed and the hidden processes. In addition to the model error of the approximate model being mitigated and correlation between the observations and hidden dimensions of the approximate model being reserved, the perfect model structure is utilized to calculate the $\cB(\u)$ in FDT with preconditioner. Therefore, the performance of the FDT with preconditioner outperforms concatenate FDT and imperfect FDT.

\section{Discussion and Conclusions}\label{Conclusion}
In this paper, the skill of a rich class of nonlinear stochastic models, known as the ``conditional Gaussian nonlinear system'' (CGNS) \eqref{CGNS}, as both a cheap surrogate model and a fast preconditioner is explored to advance many computationally challenging tasks in complex nonlinear systems.
The CGNS not only preserves the main underlying physics of nature but reproduces the observed intermittency, extreme events and other non-Gaussian features as well. The closed analytic formula of solving the conditional statistics facilitates the development of many mathematical theories and fast numerical algorithms. Three topics are covered in this paper. First, the closed analytic formulae of the conditional statistics of the CGNS allow an efficient and accurate data assimilation scheme. It is shown in Section \ref{Sec:DA_Forecast} that the data assimilation skill of a suitable CGNS approximate forecast model outweighs the EnKBF applying directly even to the perfect model in the presence of strong nonlinear and turbulent features. The latter may suffer from filter divergence when the observational process is highly nonlinear with small observational uncertainties. Second, as is shown in Section \ref{Sec:Parameter_Estimation}, the CGNS allows the development of a fast algorithm for simultaneously estimating the parameters and the unobserved state variables with uncertainty quantification in the presence of only partial observations. Utilizing an appropriate CGNS as a preconditioner significantly reduces the computational cost in accurately estimating the parameters in the original complex system. The same CGNS can also serve as a surrogate model for reproducing the large-scale dynamics and statistics of nature and can be applicable to ensemble forecast. Finally, the CGNS advances rapid and statistically accurate algorithms for both computing the probability density function and sampling the trajectories of the unobserved state variables. Showing in Section \ref{Sec:Response}, these fast algorithms facilitate the development of an efficient and accurate data-driven method for predicting the linear response of the original system with respect to parameter perturbations based on a suitable CGNS preconditioner.

Several important topics are remained as future work. First, it is crucial to systematically determine the CGNS approximate model. A promising approach is to write down the general structure of \eqref{Split_Form_CGNS} and then apply a parameter estimation algorithm together with certain sparse identification method to prevent the overfitting issue. Second, it is of practical significance to develop further pathways that allow to apply the CGNS to more complicated systems and explore the approximate errors. Incorporating the CGNS into intermediate complicated models or even certain versions of the general circulation models (GCMs), say for climate science or geophysics, can be interesting and practically useful tasks. Finally, the CGNS may have the potential to combine with machine learning algorithms to advance the ensemble forecast. One possible direction is to exploit the analytic formula in \eqref{Joint_X_Y} for the ensemble forecast, where the complicated nonlinear interactions in solving the conditional distributions can be replaced by a cheaper machine learning architecture.

\section*{Acknowledgments}
The research of N.C. is partially funded by the Office of VCRGE at UW-Madison and ONR N00014-21-1-2904. Y. L. is supported as a graduate student under the these grants. The research of Y.L. was also supported in part by NSF Award DMS-2023239 through the IFDS at UW-Madison. The work of H.L. is partially funded by NSF Award DMS-2108856.

\appendix

\section{Details of the EM algorithm for parameters estimation}\label{sec: SI_EM}
\subsection{The EM algorithm for the CGNS}
In this appendix, we provide some technical details about Algorithm~\ref{alg: EM_CGNS}, which concerns the use of CGNS as a preconditioner in the EM approach for parameters estimation. To fix ideas, we will focus on a special case in which the original nonlinear system \eqref{Split_Form_General_Eqn} is subject to diagonal and additive noise, which is sufficient for the applications considered in Section~\ref{Sec:Parameter_Estimation}. For the simplicity of discussions, we will also consider only the real-valued variables and systems. We refer to \cite{chen2020learning} for more general settings; see also \cite{ghahramani1999learning}.

We first generate the discrete partial observations $\widehat{\X}$ using for instance the Euler-Maruyama scheme~\cite{gardiner2004handbook} applied to the original nonlinear system~\eqref{Split_Form_General_Eqn}, with a sufficiently small time step size $\Delta{t}$. Given a CGNS approximation of the form \eqref{CGNS}, we explain now how the conditional distribution $p^{M}(\widehat{\Y} | \widehat{\X}, \btheta^M_{k - 1})$ and the cost function $\widetilde{\cQ}(\btheta^M;\btheta^M_{k - 1})$ in lines 5-6 of Algorithm~\ref{alg: EM_CGNS} are formed at each EM iteration for the CGNS.

Due to the above choice of the noise terms in \eqref{Split_Form_General_Eqn}, $\B_{1}$ and $\b_{2}$ in its corresponding CGNS approximation \eqref{CGNS} are simply diagonal matrices with unknown diffusion coefficient parameters appearing on the corresponding main diagonal. In the following, we also denote by $\bxi$ the parameters in the drift part of the CGNS. Then, the collection $\btheta^M$ of the parameters in the CGNS is $\btheta^M = (\bxi, \mathrm{diag}(\B_{1}), \mathrm{diag}(\b_{2}))^\top$. The discretization of \eqref{CGNS} using the Euler-Maruyama scheme reads,
\begin{subequations}\label{eq: discrete forward}
    \begin{align}
    \X^{j+1} = & \X^j + (\A_0(\X^j, t; \bxi) + \A_1(\X^j, t; \bxi) \Y^j) \Delta t + \B_{1}  \sqrt{\Delta t} \bm{\varepsilon}^j_{1}, \\
    \Y^{j+1} = & \Y^j + (\a_0(\X^j, t; \bxi) + \a_1(\X^j, t; \bxi)\Y^j) \Delta t + \b_{2} \sqrt{\Delta t} \bm{\varepsilon}^j_{2},
    \end{align}
\end{subequations}
where $j = 0, \ldots, J$ for some fixed positive integer $J$. Here, $\bm{\varepsilon}^j_{1}$ and $\bm{\varepsilon}^j_{2}$ are independent and identically distributed Gaussian white noises. They have the same dimensions as $\X$ and $\Y$, respectively, due again to the way the noise forcing is chosen. Assume also all the parameters in the drift part appear as multiplicative prefactors of some functions of $\X^j$ and $\Y^j$.

We introduce now some additional notations in order to put \eqref{eq: discrete forward} into a more compact form to be used below. Denote by $\M^j$ the matrix that includes those linear/nonlinear functions in the drift part which are multiplied by the parameters $\bxi$. Denote also by $\S^j$ those terms that do not involve parameters such as the first terms $\X^j$ or $\Y^j$ in~\eqref{eq: discrete forward}. Finally, let $\R$ be the covariance matrix associated with the noise terms in \eqref{eq: discrete forward}; namely, the diagonal matrix whose diagonal consisting of those from the diagonal matrices $\B_1\B^\top_1 \Delta t$ and $\b_2\b^\top_2\Delta t$. Apparently, there is a one-to-one correspondence between the diagonal of $\R$ and the parameters in the diffusion terms. With these notations we can rewrite \eqref{eq: discrete forward} into
\begin{equation}\label{eq: discrete forward_v2}
    \u^{j+1} = \M^j \bxi  + \S^j + \R^{1/2} \bm{\varepsilon}^j, \quad j = 0, \ldots, J,
\end{equation}
where $\u^{j+1} = (\X^{j+1}, \Y^{j+1})^\top$ and $\bm{\varepsilon}^j = (\bm{\varepsilon}_1^j, \bm{\varepsilon}_2^j)^\top$. Thus, at each time step, given $\M^j$ and $\S^j$, $\u^{j+1}$ follows a Gaussian distribution with mean $\bmu^j = \M^j \bxi  + \S^j$ and variance given by $\R$:
\begin{equation}
p^M(\u^{j+1} | \M^j, \S^j, \btheta^M) = \frac{1}{\sqrt{(2\pi)^N}} |\R|^{-\frac{1}{2}} \exp \Big(- \frac{1}{2} (\u^{j+1} - \bmu^j)^\top (\R)^{-1} (\u^{j+1} - \bmu^j)\Big),
\end{equation}
where $N$ is the dimension of the phase space.

At the $k$-th EM iteration for each $k=1,\ldots, K$, the parameters $\btheta^M_{k-1}$ is already computed. We can thus use \eqref{Smoother_Main} to compute the optimal smoother estimate $p^M(\mathbf{Y}(t)|\mathbf{X}(s), s\in[0,T], \btheta^M_{k-1})$ or equivalently its ``discretized'' form $p^M(\Y^j|\widehat{\X}, \btheta^M_{k-1})$ for $j = 0,\ldots, J$ In the E-Step. 
On the other hand, by exploiting the relationship in \eqref{eq: discrete forward_v2} for all $j$, the M-Step (line 6 in  Algorithm~\ref{alg: EM_CGNS}) is solved via minimizing the following cost function:
\begin{equation}\label{eq: SI_obj}
  \widetilde{\cL} =  \frac{1}{2} \sum_{j=J_1}^{J-1} \left\langle (\u^{j+1} - \M^j \bxi - \S^j)^\top (\R)^{-1} (\u^{j+1} - \M^j \bxi - \S^j) \right\rangle + \frac{J'}{2} \log |\R|,
\end{equation}
where the summation of $j$ starts from a certain non-zero integer $J_1$ to eliminate the inaccuracy from the burn-in period and $J' = J-J_1$.  Note that \eqref{eq: SI_obj} corresponds to the negative of $\cQ(\btheta; \btheta_k)$ defined by \eqref{eq: E-step} if $J_1 = 1$, after dropping some constant terms independent of $\bxi$ and $\R$.
In~\eqref{eq: SI_obj}, $\left\langle \cdot \right\rangle$ denotes the expectation over the uncertain component of $\u^j$ and $\u^{j+1}$, namely $\Y^j$ and $\Y^{j+1}$ while the expectations of the observed component $\X^j$ and $\X^{j+1}$ are simply themselves since the time series of them are given.  Since the hidden variables $\Y^j$ appears in a linear way in the matrix $\M^j$ due to the structure of the CGNS, only the quadratic terms of $\Y^j$, namely $\langle \Y^{j + 1}, (\Y^{j + 1})^\top \rangle$, $\langle \Y^{j + 1}, (\Y^{j})^\top \rangle$, $\langle \Y^{j}, (\Y^{j})^\top \rangle$, need to be solved in the expectation in~\eqref{eq: SI_obj}. These terms can be solved via some manipulations of the results from the closed formulae of the smoother estimates \eqref{Smoother_Main}. Details can be found in~\cite[Appendix~A]{chen2020learning}.
To find the minimum of~$\widetilde{\cL}$, we set $\frac{\partial \widetilde{\cL}}{\partial \xi_i} =0$ and $\frac{\partial \widetilde{\cL}}{\partial R_{\ell\ell}} =0$ for each component $\xi_i$ of $\bxi$ and each diagonal element $R_{\ell\ell}$ of $\R$, which leads to
\begin{subequations}\label{eq: SI_Equation_R_Theta}
    \begin{align}
        \R &= \frac{1}{J'} \sum_{j=J_1}^{J-1} \left\langle(\u^{j+1} - \M^j \bxi - \S^j)(\u^{j+1} - \M^j\bxi - \S^j)^\top\right\rangle, \\
        \bxi &= \D^{-1} \c,
    \end{align}
\end{subequations}
where
\begin{equation}\label{eq: SI_aux_physics_constraint}
	\D = \sum_{j=J_1}^{J-1} \left\langle(\M^j)^\top\mathbf{R}^{-1}\M^j\right\rangle \quad \text{and} \quad \c = \sum_{j=J_1}^{J-1}\left\langle(\M^j)^\top\mathbf{R}^{-1}(\u^{j+1} - \S^j)\right\rangle.
\end{equation}
Note that we solve~\eqref{eq: SI_Equation_R_Theta} based on an iteration method where~$\bxi$ is obtained given $\R$ from the previous step, and then $\R$ is calculated from the updated $\bxi$.

\subsection{The last M-Step in Algorithm~\ref{Alg:EM}}
Recall Algorithm \ref{alg: EM_CGNS}. The last M-Step (line 8 in the algorithm) requires to solve the minimization of the cost function, which is based on the original nonlinear system. The main difference here compared with minimizing the cost function associated with the CGNS (line 6 in the algorithm) is that higher order moments of $\Y$, resulting from the general nonlinear structure of the original complex system, may be involved. If the nonlinearity of the original complex system is up to quadratic, then the expectation in the analogue of \eqref{eq: SI_obj} for the original system may involve up to the fourth moments of~$\Y$.

Note that these moments are calculated based on the smoother estimate from the E-Step (line 7), which still utilizes the CGNS. In other words, the smoother estimate only provides a conditional Gaussian distribution. The higher order moments are thus computed based on the quasi-Gaussian closure approximation, which are represented by the known information from the mean and the variance of the conditional Gaussian smoother estimate. For example, denote by $Y_i$ a scalar component of $\Y$. So does $Y_j, Y_k$ and $Y_m$. Then the third order moment $\langle Y_i Y_j Y_k \rangle$ and the fourth order moment $\langle Y_i Y_j Y_k Y_m \rangle$ can be obtained as follows
\begin{equation} \label{Eq_higher_moments}
  \begin{aligned}
 \langle Y_i Y_j Y_k \rangle &= \mu_i\mu_j\mu_k + \mu_k \sigma_{ij} + \mu_i \sigma_{kj} + \mu_j \sigma_{ik}, \\
 \langle Y_i Y_j Y_k Y_m \rangle &= \langle Y_i Y_j Y_k \rangle\mu_m + \mu_i\mu_j \sigma_{km} + \mu_i\mu_k \sigma_{jm} + \mu_k\mu_j \sigma_{im} + \sigma_{ij}\sigma_{km} + \sigma_{ik}\sigma_{jm} + \sigma_{jk}\sigma_{im},
  \end{aligned}
\end{equation}
where $\mu_i$ and $\sigma_{ij}$ are mean and covariance of corresponding components.

\section{Calculating $\cB(\u)$ in the linear response operator}\label{sec: Bu_form}
In light of~\eqref{eq: FDT_Formula} and~\eqref{Joint_X_Y_Equilibrium_Imperfect}, one has the following explicit expression of $\cB(u)$
\begin{equation}
\begin{aligned}
	\cB(\u) &= - \frac{\textrm{div}_{\u}(\w(\u) \peq^{M | \obs}(\u))}{\peq^{M | \obs}(\u)}\\
&= - \sum_{i = 1}^N \frac{\partial}{\partial{\u_i}}\w_i(\u) - \sum_{i=1}^N \w_i \frac{\partial}{\partial \u_i} \peq^{M | \obs}(\u)\,.
\end{aligned}
\end{equation}
When perturbing the parameters of forcing in the observed processes, the forms of $\cB (\u)$ from the perfect model and the approximate model are the same, since the parameters $F_1$ and $F_2$ appear exactly the same way as in both the perfect and the approximate model. The $\cB(\u)$ term reads as follows
\begin{equation}
	\cB(\u) = - \frac{\partial}{\partial x_1} \peq^{M | \obs}(\u) - \frac{\partial}{\partial x_2} \peq^{M | \obs}(\u)\,.
\end{equation}
When perturbing the parameters in linear interactions terms that appear in both the observed and hidden processes, the formulation of $\cB(\u)$ from the perfect and approximate models are different. Given the perturbation vector $\w(\u) = (y_1, y_2, -x_1, -x_2)^{\top}$, the $\cB(\u)$ from the perfect model is as follows
\begin{equation}
	\cB(\u) = - y_1 \frac{\partial}{\partial x_1} \peq^{M | \obs}(\u) - y_2 \frac{\partial}{\partial x_2} \peq^{M | \obs}(\u) + x_1 \frac{\partial}{\partial y_1} \peq^{M | \obs}(\u) + x_2 \frac{\partial}{\partial y_2} \peq^{M | \obs}(\u)\,.
\end{equation}
However, since there is no $L_{13}$ and $L_{24}$ parameters in the hidden processes of the approximate model, the formulation of $\cB(\u)$ from the approximate model remains
\begin{equation}
	\cB(\u) = - y_1 \frac{\partial}{\partial x_1} \peq^{M | \obs}(\u) - y_2 \frac{\partial}{\partial x_2} \peq^{M | \obs}(\u).
\end{equation}

\bibliographystyle{plain}
\bibliography{references}

\end{document}